\documentclass[10pt]{amsart}

\usepackage{amssymb,amsmath,amscd,amsfonts,a4,psfrag,graphicx,
latexsym,epsfig,color}
\usepackage[active]{srcltx}
\usepackage[all]{xy}

\everymath{\displaystyle}

\newtheorem{teo}{Theorem}[section]
\newtheorem{lem}[teo]{Lemma}

\newtheorem{exa}[teo]{Example}
\newtheorem{prop}[teo]{Proposition}
\newtheorem{defi}[teo]{Definition}

\newtheorem{conj}[teo]{Conjecture}
\newtheorem{remark}[teo]{Remark}

\newcommand{\cvd}{\hfill$\Box$}

\newcommand{\R}{\mathbb{R}}

\newcommand{\Z}{\mathbb{Z}}

\newcommand{\Aa}{{\mathcal A}}
\newcommand{\Bb}{{\mathcal B}}

\newcommand{\Ff}{{\mathcal F}}

\newcommand{\Hh}{{\mathcal H}}

\newcommand{\Mm}{{\mathcal M}}
\newcommand{\Nn}{{\mathcal N}}

\newcommand{\Pp}{{\mathcal P}}

\newcommand{\Ss}{{\mathcal S}}
\newcommand{\Tt}{{\mathcal T}}

\newcommand{\Vv}{{\mathcal V}}

\newcommand{\PP}{{\bf P}}

\newcommand{\TG}{{\mathfrak T}}
\newcommand{\hG}{{\mathfrak h}}
\newcommand{\vG}{{\mathfrak v}}
\newcommand{\iG}{{\mathfrak i}}
\newcommand{\SG}{{\mathfrak S}}
\newcommand{\pG}{{\mathfrak p}}
\newcommand{\mG}{{\mathfrak m}}

\newcommand{\Dim}{{\it Proof.\ }}

\def\cvd{\hfill$\Box$}

\title{On ideal triangulations of surfaces up to branched transit equivalences}

\begin{document}
\author{Riccardo Benedetti}
\address{Dipartimento di Matematica, Largo Bruno Pontecorvo 5, 56127 Pisa, Italy} 
\email{riccardo.benedetti@unipi.it}
\date{\today}

\begin{abstract}
We consider triangulations of closed surfaces $S$ with a given set of
vertices $V$; every triangulation can be branched that is enhanced to
a $\Delta$-complex. Branched triangulations are considered up to the
$b$-transit equivalence generated by $b$-flips (i.e. branched diagonal
exchanges) and isotopy keeping $V$ pointwise fixed. We extend a well
known connectivity result for `naked' triangulations; in particular in
the generic case when $\chi(S)<0$, we show that branched
triangulations are equivalent to each other if $\chi(S)$ is even,
while this holds also for odd $\chi(S)$ possibly after the complete
inversion of one of the two branchings. Moreover we show that under a
mild assumption, two branchings on a same triangulation are connected
via a sequence a inversions of ambiguous edges (and possibly the total
inversion of one of them). A natural organization of the $b$-flips in
subfamilies gives rise to restricted transit equivalences with non
trivial (even infinite) quotient sets. We analyze them in terms of
certain preserved structures of differential topological nature
carried by any branched triangulations; in particular a pair of
transverse foliations with determined singular sets contained in $V$,
including as particular cases the configuration of the vertical and
horizontal foliations of the square of an Abelian differential on a
Riemann surface.
\end{abstract}

\maketitle

\section{Introduction} 
Different notions of ``decorated ideal triangulations'' of
$3$-manifolds (branched, pre-branched, weakly branchedp)
considered up to various transit equivalences naturally arise in the
developments of quantum hyperbolic geometry (see for instance 
\cite{NA}, \cite{QT},
\cite{AGT}) and in several other instances of 3D quantum
invariants based on state sums over triangulations (see \cite{NA}).
To understand the intrinsic content of the corresponding quotient sets
is an interesting and non trivial task.  This note arises as a simpler
but non obvious 2D counterpart of similar questions, which also
emerged within \cite{NA}, Section 5,  in a so called ``holographic" approach to 3D non
ambiguous structurers.

\medskip

Let $(S,V)$ be a compact closed connected smooth surface $S$ with a
set $V$ of $n$ marked points and Euler-Poincar\'e characteristic $\chi(S)$, such that
$\chi(S)- n < 0 $.  It is well known that $(S,V)$ carries {\it
  ideal triangulations}, say $T$.  This means that $T$ is a possibly
{\it loose} triangulation (self and multiple edge adjiacency being
allowed) whose set of vertices coincides with $V$. Clearly such ideal
triangulations of $(S,V)$ share the same numbers of edges and
triangles, $3(n-\chi(S))$ and $2(n-\chi(S))$ respectively.  It is
sometimes useful to consider an ideal triangulation $T$ as a way to
realize $(S,V)$ by assembling $2(n-\chi(S))$ ``abstract'' triangles by
gluing their ``abstract'' edges in pairs in such a way that no edge
remains unglued.  Ideal triangulations of $(S,V)$ are considered up to
the {\it ideal transit equivalence} which is generated by isotopy
fixing $V$ pointwise
and the elementary {\it diagonal exchange move} also called {\it
  flip}. Denote by $\Tt^{id}(S,V)$ the corresponding quotient set.
The following is an important well known  connectivity result.

\begin{teo}\label{naked} 
For every $(S,V)$, $\Tt^{id}(S,V)$ consists of one point.
\end{teo}

Proofs are available in several papers such as \cite{Hat}, \cite{La},
\cite{Mos}, \cite{T-W}.
\smallskip

Every ``naked'' triangulation $T$ carries some {\it branchings}
$(T,b)$ (see Lemma \ref{b-exist}), where by definition $b$ is a system
of edge orientations which lifted to every abstract triangle $(t,b)$
of $T$ is induced by a (local) ordering of the vertices, so that every
edge  goes towards the bigest endpoint; equivalently $b$ promotes
$T$ to be a $\Delta$-{\it complex} accordingly with \cite{HATCHER},
Chapter 2. It is easy to see that for every branching $(T,b)$, every
naked flip $T\to T'$ can be enhanced to some {\it $b$-flip} $(T,b)\to
(T,b')$ such that every ``persistent'' edge in both $(T,b)$ and
$(T',b')$ keeps the same orientation.  Isotopy relatively to $V$ and
$b$-flips generate the so called {\it ideal $b$-transit equivalence}
and we denote by $\Bb^{id}(S,V)$ the corresponding quotient set.  We
define the {\it symmetrized relation} by adding to the generators che
{\it complete inversion} that is we stipulate that every $(T,b)$ is
equivalent to $(T,-b)$ where $-b$ is obtained by inverting all edge
orientations of $b$, and we denote by $\tilde \Bb^{id}(S,V)$ the
corresponding quotient sets. It is not hard to see that by setting
$\sigma([T,b])=[(T,-b)]$ it is well defined an involution on
$\Bb^{id}(S,V)$ and that $\tilde \Bb^{id}(S,V) \sim
\Bb^{id}(S,V)/\sigma$.  By the topological homogeneity of every
surface, the cardinality of $\Bb^{id}(S,V)$ only depends on the
topological type of $S$ and the number $n=|V|$; sometimes we will
write $(S,n)$ instead of $(S,V)$. The following branched version of
the above connectivity result is a main result of the present note.

\begin{teo}\label{main} 
(1) If $S$ is orientable or if it is non orientable and $\chi(S)$ is
  even and strictly negative, then for every $(S,V)$, $\Bb^{id}(S,V)$
  consists of one point.

(2) If $S$ is not orientable and either $\chi(S)=0$ or $\chi(S)$ is
  odd, then for every $(S,V)$, $\tilde \Bb^{id}(S,V)$ consists of one
  point.

\end{teo}

As $\tilde \Bb^{id}(S,V)$ is a quotient of $\Bb^{id}(S,V)$ by an
involution, it follows that in case (2), $| \Bb^{id}(S,V)|\leq 2$.
\begin{conj}
  If $S$ is not orientable and either $\chi(S)=0$ or $\chi(S)$ is odd,
  then for every $(S,V)$, $|\Bb^{id}(S,V)|=2$.
\end{conj}

This will be confirmed at least for $\Bb^{id}(\PP^2(\R),2)$ (Proposition \ref{PP-2}).

\smallskip

Assuming Theorem \ref{naked}, we will provide two constructive proofs
of Theorem \ref{main} each with its subtleties and constructions of
distinguished $b$-transits.  A key ingredient will the move of {\it
  inverting an ambiguous edge} (see Section \ref{inversion-amb} and in
particular Theorem \ref{inversive}).

By adding to the $b$-flips the branched positive $0\to 2$ $b$-{\it bubble}
moves and their inverse (or equivalently the stellar $1\to 3$
branched moves and their inverse), we get the {\it completed
  $b$-transit equivalence} with quotient set denoted by $\Bb(S,V)$.  A
positive bubble move produces an ideal triangulation of $(S,V')$ where
$V'$ contains one further marked point of $S$; if it is part of a
$b$-transit which connects two ideal branched triangulations of
$(S,V)$, then it must be compesated later by a negative inverse move.
We will see a quick direct proof of the following weaker connectivity
result (no matter if $S$ is orientable or not).

\begin{prop}\label{complete} For every $(S,V)$, $\Bb(S,V)$ consists of one point.
\end{prop}

\medskip

We will see in Section \ref{generalia} that $b$-flips 
can be naturally organized in some sub-families so that
more restrictive transit equivalences can be defined, with non-trivial
(actually infinite, see Remark \ref{infinite-s}) quotient sets.  Another main theme (Sections
\ref{sliding}) is to point out the {\it intrinsic content} of
these various transit equivalences, that is some relevant
structures on $(S,V)$, carried by every $(T,b)$, which are invariant
under a given instance of transit equivalence.  Theorem \ref{main}
itself should be enlighten by the mutations of such structures along
any ideal $b$-transit.

\medskip

{\bf The dual viewpoint.} Let $S_V$ be the surface with $n$ boundary
components ($n=|V|$) obtained by removing from $S$ a small open ball
around each $v\in V$. For every ideal triangulation $T$ of $(S,V)$,
the $1$-skeleton $\theta=\theta_T$ of the dual cell decomposition is a
generic (internal) spine of $S_V$. In fact $\theta$ is a graph with
$3$-valent vertices and $S_V$ is a ribbon graph which tickens
$\theta$. If $(T,b)$ is branched, this promotes $\theta$ to be a {\it
  transversely oriented train track} $(\theta,b)$ - for simplicity we
keep the same notation ``$b$''.

\begin{remark}\label{oriented}{\rm
    If $S$ is {\it oriented}, then $(\theta, b)$ can be equivalently
    considered as an {\it oriented} train track by means of the
    following dual orientation convention:
\medskip

{\it At every transversal intersection point of $T$ and $\theta$, 
an oriented edge of $(T,b)$
followed by the dual oriented branch of $(\theta,b)$ realize 
the orientation of $S_V$ (that is every intersection number is equal to $1$).}  
}
\end{remark}

By definition, $(\theta,b)$ is a {\it (transversely oriented) branched
  spine} of $S_V$.  Viceversa, every ribbon graph, $\bar S$ say,
carried by a (possibly branched) spine $\theta$ as above gives rise to
a (possibly branched) ideal triangulation $T=T_\theta$ of $(S,V)$
obtained by filling each boundary component of $\bar S$ with a
punctured $2$-disk.  Flips and bubbles, possibly branched, can be
equivalently rephrased in terms of (branched) spine moves. We will
freely adopt both equivalent dual viewpoints.

\begin{remark} {\rm Although they are equivalent, there is some
    qualitative difference beetween spines and triangulations. A flip
    is a {\it discrete} transition with a cell decomposition as
    intermediate ``state'' which is no longer a triangulation (it
    includes one quatrilateral). The corresponding spine transition
    can be realized by a {\it continuous} deformation passing through
    a non generic spine (with one $4$-valent vertex).}
\end{remark}

\section{Generalities on $b$-transit }\label{generalia}
An ``abstract'' $b$-flip acts on a quadrilateral $Q$ endowed with a
branched triangulation $(t_1\cup t_2, b)$ made by two triangles with
one common edge $e= t_1\cap t_2$ (a diagonal of the quadrilateral).  A
$b$-flip produces another branched triangulation $(t_1'\cup t_2',b')$
of $Q$ made by two triangles having as common edge $e'=t_1'\cap t_2'$
the other diagonal of $Q$, while $b$ and $b'$ coincide on the
persistent edges which form the boundary of $Q$.  An abstract $b$-flip
can be applied at every couple of abstract triangles of any branched
ideal triangulation $(T,b)$ of any $(S,V)$, (partially) glued in $T$
along a common edge. When we say that a $b$-flip verifies a certain
property we mean that this holds ``universally'' for every $(S,V)$ and
every triangulation $(T,b)$ at which the flip operates.
 
 \subsection{A combinatorial classification of $b$-flips}\label{flip-class}
 For every branched triangulation $(t_1\cup t_2, b)$ of $Q$ as above
 there are either one or two ways to enhance the naked flip $t_1\cup
 t_2 \to t_1'\cup t_2'$ to a $b$-flip $(t_1\cup t_2,b) \to (t_1'\cup
 t_2', b')$.  This last is sometimes denoted by $f_{e,b,b'}$ while the
 underlying naked flip is denoted by $f_e$.  Then we can distinguish a
 few families of $b$-flips.  The classification and even the
 terminology below could sound a bit arbitrary at this point. This
 will be clarified later.

 \begin{figure}[ht]
\begin{center}
 \includegraphics[width=9.5cm]{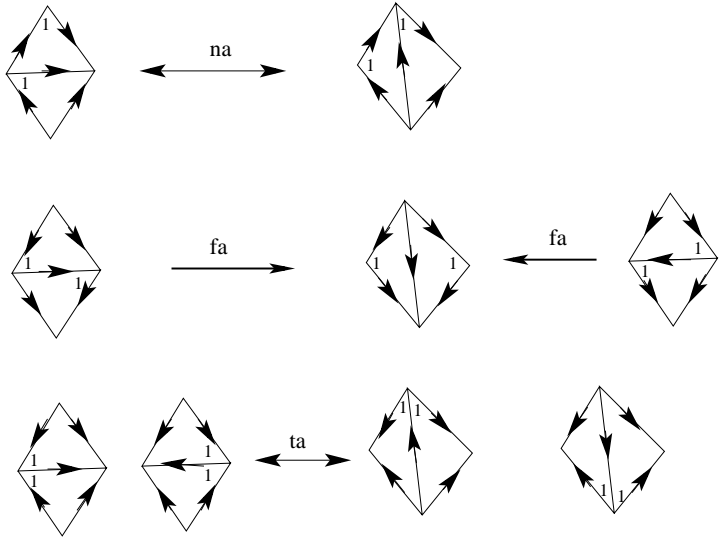}
\caption{\label{b-flip} Branched flips.} 
\end{center}
\end{figure}

  \begin{defi}\label{bflip-class}{\rm  
\begin{enumerate}
 \item   A $b$-flip $f_{e,b,b'}$
 is  {\it forced} if it is the unique
 branched flip which enhances $f_e$, starting from
 $(t_1\cup t_2,b)$.
 
 \item A $b$-flip $f_{e,b,b'}$ is {\it non ambiguous}
 if both  $f_{e,b,b'}$ and the inverse $b$-flip  $f_{e',b',b}$ are forced.
 
 \item A $b$-flip $f_{e,b,b'}$ is {\it forced ambiguous}  if it is forced
but the inverse $b$-flip is not.
 
\item A $b$-flip $f_{e,b,b'}$ is said a {\it sliding flip} ($s$-{\it flip})
if at least one among  $f_{e,b,b'}$ or its inverse  $f_{e',b',b}$ is forced.

\item A $b$-flip $f_{e,b,b'}$  is {\it totally ambiguous} (also called a
{\it bump flip}) if noone among $f_{e,b,b'}$ and $f_{e',b',b}$ is forced.

\end{enumerate}
  }
 \end{defi}

In Figure \ref{b-flip} we show typical samples of
$b$-flips in accordance with the above classification.
We have labelled by $1$ the corner of each branched triangle
 formed by the two edges that carry the {\it prevalent
  orientation}. Here $1$ is just a highlitghting label. For its meaning 
  If $1\in \Z/2\Z$, we refer to Section 5 of \cite{NA}.
 We note that the above classification of the $b$-flips is {\it
 invariant under total inversion}.

\begin{figure}[ht]
\begin{center}
 \includegraphics[width=9cm]{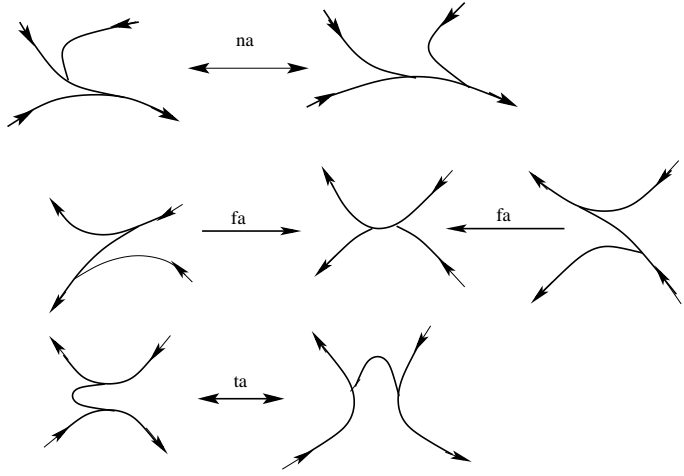}
\caption{\label{sliding_flip} Sliding and bump branched 
spine flips.} 
\end{center}
\end{figure} 

To stress it, in Figure
\ref{sliding_flip} we show the dual pictures in terms of branched
spines, provided that we have performed the total inversion on the
$b$-flips of Figure \ref{b-flip}; instead of the transverse
orientations, we prefer (as it is easier) to indicate the {\it local}
orientations on the dual train-tracks, by stipulating that all these
pictures are planar, the plane $\R^2$ is oriented by the standard
basis and we apply the orientation convention fixed in Remark
\ref{oriented}.
  
The $s$-flips as well as the $na$-flips (together with isotopy
relative to $V$) generate restricted $s$- and $na$-{\it ideal transit
  equivalence} with respective quotient sets denoted $\Ss^{id}(S,V)$
and $\Nn\Aa^{id}(S,V)$.  Clearly there are surjective projections
$$\Nn\Aa^{id}(S,V)\to \Ss^{id}(S,V)\to \Bb^{id}(S,V) $$
in particular the last quotient map is obtained by adding the bump $b$-flips
to the sliding ones.

\smallskip

{\bf Some characterizations}.
For every branched triangulation $(T,b)$ of $(S,V)$, for every
vertex $v$ of $T$, the number of corners at $v$ in its star labelled
by $1$ (as above) is even, say $2d_b(v)$. In fact given the star of $v$ an auxiliary
orientation, the $1$-labelled cornes at $v$ belong to triangles $(t,b)$
whose $b$-orientations alternate with respect to the reference one.
It is clear that
$$ \chi(S) = |V|- \sum_v d_b(v) = \sum_v (1-d_b(v)) \ . $$
We easily have that $b$-flips are characterized by the following  property.
\begin{prop}\label{s-d(v)}
 A (abstract)  $b$-flip $f_{e,b}$ is an $s$-flip if and only if for every
  $(S,V)$ and for every application of the flip on triangulations of
  $(S,V)$, $(T,b)\to (T',b')$, we have that for every $v\in V$,
  $d_b(v)=d_{b'}(v)$.
\end{prop}

Note that when $|V|=1$ the conclusion 
 holds for {\it every} $b$-flip, not necessarily
an $s$-flip, but this is not  ``universally'' true.

Suppose now that $S$ is {\it oriented}. The orientation
corresponds the a unique simplicial {\it fundamental $\Z$-$2$- cycle}
$$ f(T,b)=\sum_t  *_{(t,b)} (t,b) $$
where every $*_{(t,b)} \in \{\pm1\}$.
Set 
$$\epsilon_\pm= \epsilon_\pm (T,b):=|\{t ; \ *_{(t,b)}= \pm 1\}| \ . $$ 
We have
\begin{lem} For every $(T,b)$, $\epsilon_+ = \epsilon_-  $.
\end{lem}
\Dim Around every $v\in V$, the $2d(v)$ $1$-colored corners 
necessarily belong to triangles with alternating signs. As every
triangle contains one $1$-colored corner the result follows.

\cvd
\smallskip

The above Lemma corroborates the validity of
Theorem \ref{main}. 
For we easily check that every $b$-flip preserves the value of
$\epsilon_+ - \epsilon_- $.  On the other hand
$$(\epsilon_+ - \epsilon_-)(T,b)= -(\epsilon_+ - \epsilon_-)(T,-b) \ . $$
\smallskip

Denote by $S_\pm= S_\pm(T,b)$ the
union of triangles such that $*_{(t,b)}= \pm 1$.
Then $S$ decomposes as $$S=S_+ \cup S_- \ . $$
Denote by $\partial S_\pm $ the boundary $1$-cycle
of the simplicial $2$-cochain supported  by $S_\pm$.

We have

\begin{prop}\label{S+}   
  A $b$-flip $f_{e,b}$ is non ambiguous if and only if it is
  an $s$-flip and for every $(S,V)$ and for every application of the
  flip on triangulations of $(S,V)$, $(T,b)\to (T',b')$, we have that
  $S_\pm(T,b)=S_\pm(T',b')$, hence also $\partial S_\pm$ 
  is preserved.
\end{prop}

\medskip 

\begin{remark}\label{on-bubble} {\rm Also branched bubble and stellar $1\leftrightarrow 3$ moves
admit a {\it sliding vs bump} classification. This leads to the corresponding 
completed sliding equivalence with quotient set $\Ss(M)$ which projects onto $\Bb(M)$ (see \cite{NA}, Section 5.3).}
\end{remark}

\subsection{Inversion of an ambiguous edge}\label{inversion-amb}

\begin{defi}\label{amb-edge}{\rm 

(1) Let $T$ be a naked ideal triangulation of $(S,V)$.  An edge $e$ of $T$ is
      said {\it trapped} if it results by the identification of two
      edges of {\it one} ``abstract'' triangle. Otherwise, $e$ is said
      {\it untrapped}. A trapped edge corresponds to a one vertex loop
      in the dual spine $\theta$.

(2) Given a branching $(T,b)$, a $b$-oriented
edge $e$ is said {\it ambiguous in $(T,b)$} if by inverting its
orientation we keep a branched triangulation $(T,b')$. 
}
\end{defi}

\begin{lem}\label{eliminate-amb} 
  If $e$ is ambiguous and untrapped in $(T,b)$, then $(T,b)$ and
  $(T,b')$ (as in (2) of Definition \ref{amb-edge}) are connected by two
  $b$-flips.
\end{lem}
\Dim Denote by $f_{e,b,b"}$ a $b$-flip that enhances the naked flip
$f_e$ with inverse naked flip $f_{e'}$. Then we easily see that the
untrapped edge $e$ is ambiguous if and only if either $f_{e,b,b"}$ is
forced ambiguous or it is totally ambiguous.  Hence $f_{e,b,b"}$
followed by $f_{e',b",b'}$ convert $(T,b)$ to $(T,b')$.

\cvd

\smallskip

Then {\it we can add the elementary move of inverting any untrapped
    ambiguous edge without changing the ideal $b$-transit
    equivalence.}
     
\begin{remark}\label{trapped} {\rm  Every trapped edge $e$ of $(T,b)$ is ambiguous but
    there is not any apparent {\it local} sequence of b-flips that
    inverts $e$.}
\end{remark}

\smallskip

We have

\begin{lem} \label{no-trapped} (1) For every $T$ as above there is a
  sequence of flips $T\Rightarrow T'$ such that $T'$ does not contain
  trapped edges.  (2) If $T$ and $T'$ do not contain trapped edges
  then they can be connected by a sequence of flips through
  triangulations without trapped edges.
  \end{lem} 
  
 \Dim The vertex of a loop in the spine $\theta$ which is dual to a
 trapped edge of $T$ is connected by an edge to the rest of the spine.
 By performing the dual flip at this edge we remove the loop without
 introducing new ones. If such a loop appear in a sequence of flips
 connecting $T$ and $T'$ as in (2), then we can follow it till it
 disappears so thet we can eventually remove it from the sequence.

\cvd

\subsection{Bad nutshells}\label{bad-nutcshell}
\medskip

A positive naked $0\to 2$ bubble produces a so called {\it nutshell}
made by two triangles identified along two common edges.
Not every branched nutshell $(N,b)$ supports a negative $2\to 0$ $b$-bubble.

\begin{defi}\label{nutshell} {\rm A branched nutshell
$(N,b)$ is {\it bad} if the two boundary edges form an {\it oriented}
circle. Otherwise $(N,b)$ is a {\it good} nutshell.}
\end{defi}

\medskip

The following Lemma is immediate.

\begin{lem}\label{2_nutshell} 
(1) If $(N,b)$ is a bad nutshell, then the central vertex
is necessarily either a pit or a source.

(2) $(N,b)$ is good if and only if
it supports a negative $b$-bubble move.

(3) Two different good nutshells $(N,b)$ and $(N,b')$
sharing the same oriented boundary edges are connected  by either
one or two consecutive inversions of internal (hence untrapped) ambiguous edges.
\end{lem}

A positive naked $1\to 3$ move produces a so called {\it triangular star}.
Similarly as before, not every branched triangular star, say $(\SG,b)$,
supports a negative $b$-$(3\to 1)$ move.

\begin{defi}\label{t-star} {\rm A branched triangular star $(\SG,b)$
is {\it bad} if the three boundary edges form an {\it oriented}
circle. Otherwise $(\SG,b)$ is  {\it good}.}
\end{defi}

Easily we have

\begin{lem}\label{2t-star}
  (1)  If $(\SG,b)$ is bad, then the central vertex
  is necessarily either a pit or a source.

  (2) $(\SG,b)$ is good if and only if
  it supports a negative $b$-$(3\to 1)$  move.

(3) Two good $b$-triangular stars $(\SG,b)$ and $(\SG,b')$
  sharing the same oriented boundary edges are connected  by
a finite sequence of consecutive inversions of internal (hence untrapped)
ambiguous edges.
\end{lem}

\subsection{Existence of branched triangulations}

We have

\begin{figure}[ht]
\begin{center}
 \includegraphics[width=7cm]{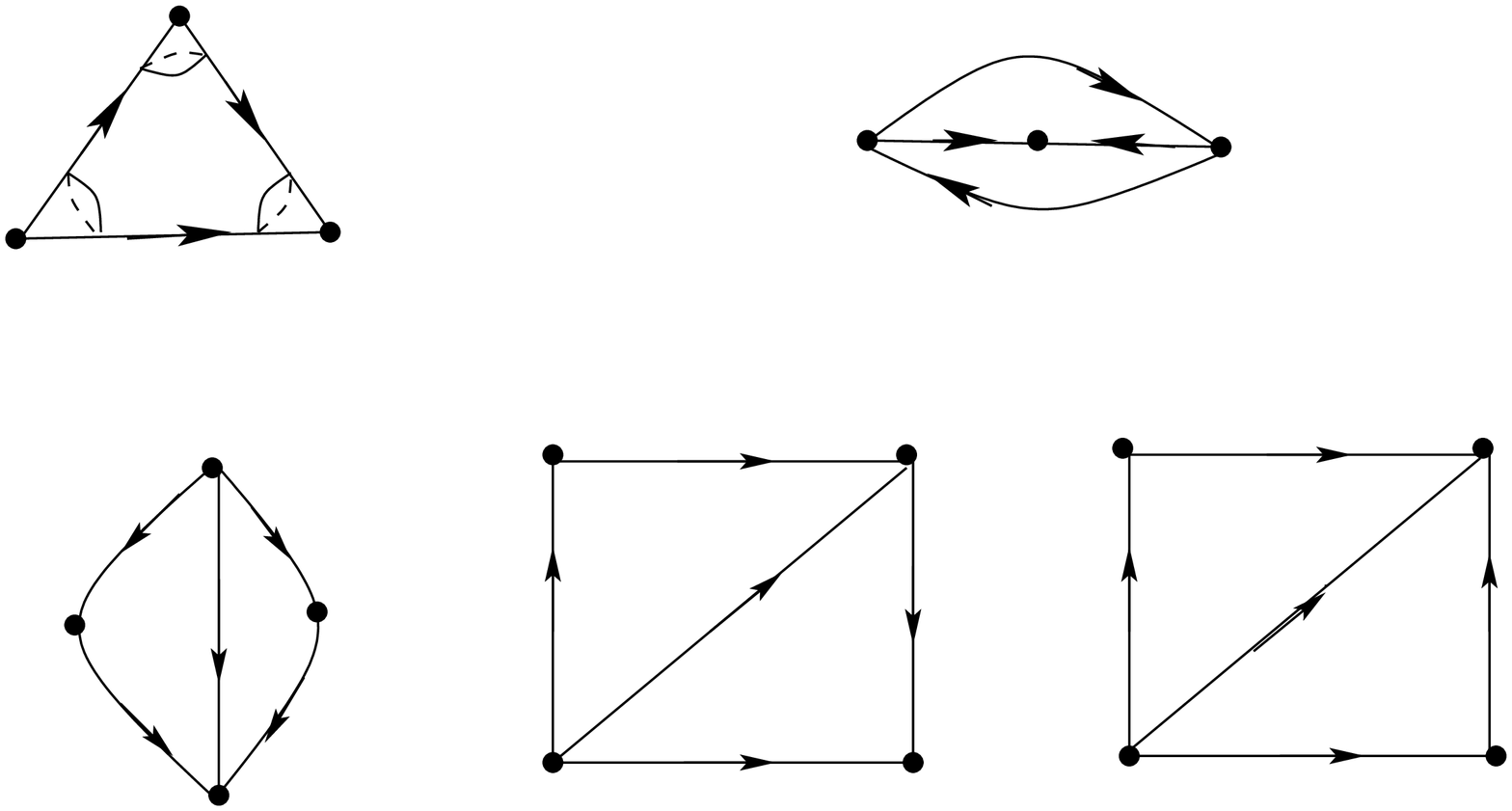}
\caption{\label{b-esiste} Existence of branched triangulations.} 
\end{center}
\end{figure} 

\begin{figure}[ht]
\begin{center}
 \includegraphics[width=7cm]{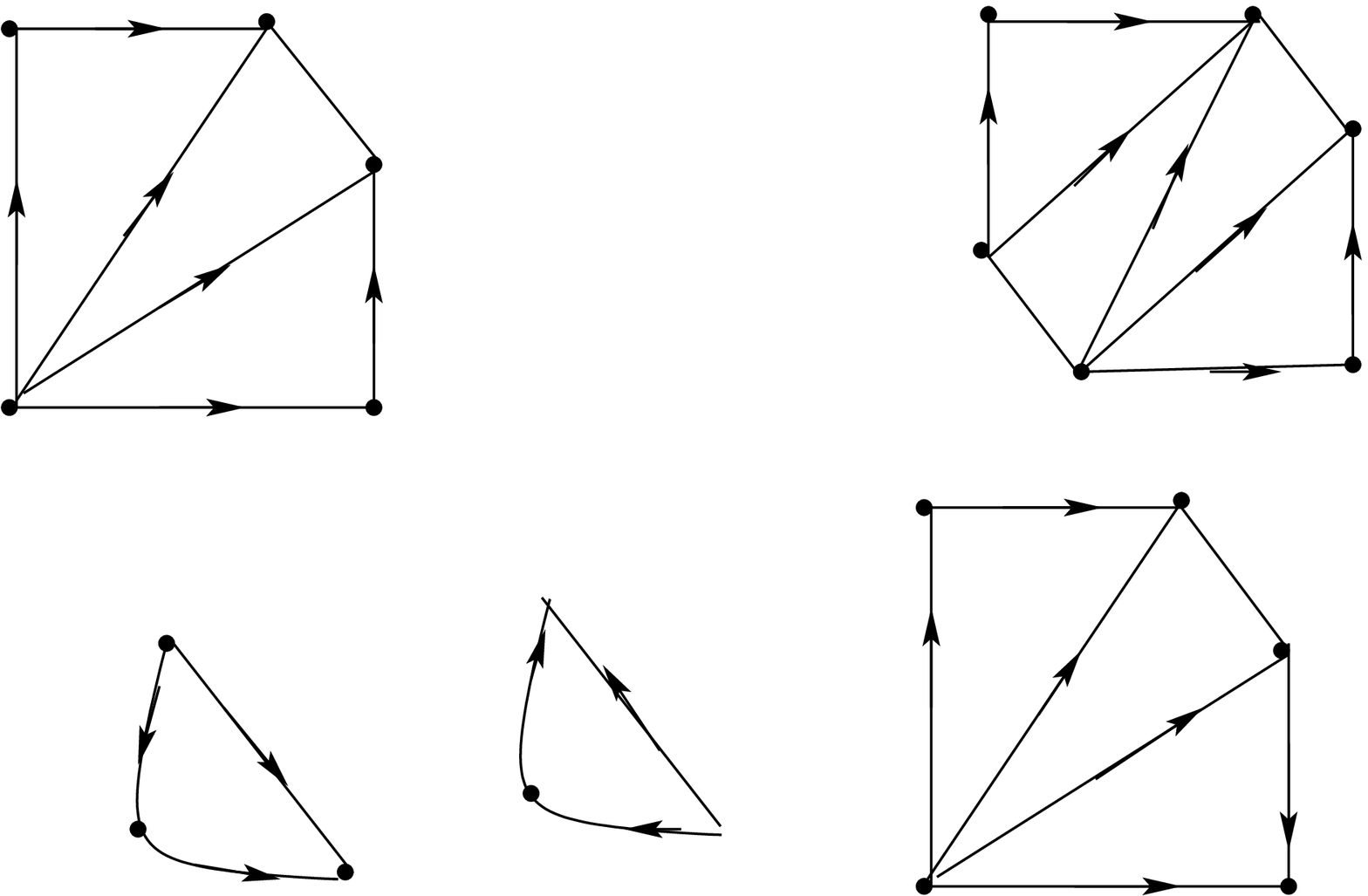}
\caption{\label{b-esiste2} Existence of branched triangulations.} 
\end{center}
\end{figure} 

\begin{lem} \label{b-exist} Every ideal triangulation $T$ of $(S,V)$ can be
  branched.
\end{lem}
\Dim Thanks to Theorem \ref{naked} and due to the fact that for every
branched triangulation $(T,b)$, every naked flip $T\to T'$ can be
enanched to a $b$-flip $(T,b)\to (T',b')$, it is enough to show that
every $(S,V)$ admits a branched triangulation. By using the
$b$-bubbles we see that if $(S,V)$ admits such a triangulation, then
this holds for every $(S,V')$ such that $|V'|\geq |V|$. Then it is
enough to show that for every $S$ there exists a branched
triangulation $(T,b)$ of $(S,V)$ such that $n=|V|$ is the minimum for
which $\chi(S)-n <0$. There are of course several ways to do it.  We
indicate a way which will be suited for further use.  If $S$ is a
sphere then $n=3$, and a branched triangulation is obtained by
gluing two copies of a branched triangle along the common boundary
(see the left side of the first row of Figure \ref{b-esiste}).  If
$S=\PP^2(\R)$ is a projective plane, then $n=2$ and we can use the
realization of $\PP^2(\R)$ by identifying the two edges of a bigon and
triangulate it with one internal vertex - see the right side of the
first row of Figure \ref{b-esiste}.  For all other $S$, $n=1$.  If
$S= \PP^2(\R)\# \PP^2$ is a Klein bottle, that is the connected sum of
two projective planes, we get a branched triangulation of $(S,V)$
starting from a realization of $S$ by identifying the boundary of a
quadrilateral obtained by gluing two ``truncated bigons" - see the
left side of the second row of Figure \ref{b-esiste}, at the middle of
the row we show another branched triangulation starting from a
realization of the bottle by identifying in pairs the opposite edges
of a quadrilateral; similarly on the right side we suggest a
triangulation of the torus $S^1\times S^1$.  In Figure \ref{b-esiste2}
we show the elementary bricks in order to realize all other cases.
These bricks are branched triangulations of certain surfaces with
boundary.  In the first row we see either a one or twice-pierced
torus that is a torus from which one has removed respectively one or
two open $2$-disks; they are obtained by means of a one or
twice-truncated quadrilateral with opposite edge identified in pairs.
In the second row, right side we see a similar realization of a
one-pierced Klein bottle; on the left a one-pierced projective plane
given by means of a truncated bigon with the two possible
branchings. The non oriented edges are ambiguous so that their
orientation can be chosen arbitrarily. If $S$ is orientable of genus
$g>1$ we can realize it by means of a chain of $g-2$ twice-pierced
tori capped by two one-pierced ones. If $S$ is non orientable and
$\chi(S)=2-r <0$ is odd, set $g=\frac{r-1}{2}$; then we can obtain $S$
by means of chain of $g-2$ twice-pierced tori capped by a one-pierced
torus and a one-pierced projective plane. If $\chi(S)-r<0$ is even,
set $g=\frac{r-2}{2}$, then we obtain $S$ by means of chain of $g-2$
twice-pierced tori capped by a one-pierced torus and a one-pierced
Klein bottle.

\cvd

\begin{defi}\label{connection-edge}{\rm The boundary of every
    brick as above is union of loops with one vertex. The corresponding
    edge of the triangulation is called a {\it connection edge}.}
 \end{defi}
    
    \smallskip
    
    For every triangulation of $S$ obtained so far, the connection edges
    become separating loops that
decompose $S$ by the bricks; in Figure \ref{b-esiste2}
these edges correspond to the ambiguous edges on the boundary of the
truncated quadrilaterals or to the non ambiguous edge in the branched
truncated bigons.
 
 \subsection{Preliminary reductions to face  the ideal 
$b$-transit equivalence}\label{reduction}
We assume Theorem \ref{naked}. Isotopy relative to $V$ will be understood.
At the end of this section we will obtain a quick proof of Proposition \ref{complete}.

\begin{lem}\label{sameT}  The following facts are equivalent
to each other:
\begin{enumerate}

\item $\Bb^{id}(S,V)$ consists of one point.

\item For every naked ideal triangulation $T$ of $(S,V)$, every two
  branchings $(T,b)$ and $(T,b')$ are connected by a chain of
  $b$-flips.

\item There exists a naked ideal triangulation $T$ of $(S,V)$ such
  that every two branchings $(T,b)$ and $(T,b')$ are connected by a
  chain of $b$-flips.
\end{enumerate}
\end{lem}

\Dim Obviously $(1)\Rightarrow (2)\Rightarrow (3)$.  In order to prove
$(3)\Rightarrow (1)$, we argue similarly to the proof of
Lemma \ref{b-exist}.
Let $(T_1,b_1)$ and $(T_2,b_2)$ be ideal
triangulations of $(S,V)$. By Theorem \ref{naked}, there are naked ideal
$b$-transits $T_j \Rightarrow T$, $j=1,2$. There is no
obstruction to enhance them to sequences of $b$-flips
$(T_j,b_j)\Rightarrow (T,b'_j)$. The Lemma follows immediately.

\cvd

Similar statements hold for both the symmetrized and completed
$b$-transit equivalences.  

\begin{defi}\label{delta}{\rm Given $(T,b)$ and $(T,b')$, denote by
    $\delta(b,b')$ the set of edges of $T$ at which $b$ and $b'$ are
    opposite. An edge $e\in \delta(b,b')$ is said {\it disoriented}.}
\end{defi}

The previous considerations suggest two possible ``strategies'' in order to
prove that $\Bb^{id}(S,V)$  (or $\tilde \Bb^{id}(S,V)$ or $\Bb(S,V)$) consits of one point.

\medskip

{\bf (A)} For a given $(S,V)$ detect a distinguished
naked triangulation $T$ for which one can check directly that (3) of Lemma \ref{sameT}
holds.

\medskip
 
{\bf (B)} To point out a few procedures such that for every couple
$(T,b)$ and $(T,b')$ such that $\delta(b,b')$ is non empty, we can
apply one of them producing $b$-transits $(T,b)\Rightarrow (T',b_1)$
and $(T,b')\Rightarrow (T',b_2)$ such that
$$ |\delta(b_1,b_2)|<|\delta(b,b')| \ . $$

\begin{remark}\label{AvsB}{\rm We understand the difference between 
the two strategies if we apply them in order to realize a $b$-transit between
two {\it arbitrary} triangulations $(T_1,b_1)$ and $(T_2,b_2)$ of some $(S,V)$.
Via {\bf A}, if $T$ is the distinguished triangulation we must preliminarly
connect both $(T_1,b_1)$ and $(T_2,b_2)$ with some $(T,b)$ and $(T,b')$
respectively, by applying twice Theorem \ref{naked}.
Via {\bf B}, it is enough to connect $(T_1,b_1)$ to some
$(T,b):= (T_2,b)$ and set $(T,b')=(T_2,b_2)$. This might be relevant in terms of
computational cost.
}
\end{remark}

Strategy {\bf B}, suitably adapted, works quickly on $\Bb(S,V)$.
\smallskip

{\bf Proof of Proposition \ref{complete}:} Given $(T,b)$ and $(T,b')$
as above perform on both triangulations a bump $b$-$(1\to 3)$ move at
every triangle of $T$ in such a way that all the new oriented edges
point toward the new internal vertex. We get in this way $(T',b_1)$
and $(T',b_2)$ such that $\delta(b_1,b_2)= \delta(b,b')$.  We realize
now that every $e\in \delta(b,b')$ is untrapped and ambiguous in both
$(T',b_1)$ and $(T',b_2)$. So we conclude by several applications of
Lemma \ref{eliminate-amb} (by the way, in the present situation every
inversion of $e$ is obtained by a sequence of two bump $b$-flips).

\cvd 

We will implement both strategies and get the two promised proofs of
Theorem \ref{main}.

\section{ A-proof of the main Theorem}\label{A-main}
By implementing strategy {\bf A} we will actually obtain stronger
results. We have:

\begin{teo}\label{main-A}  (i) For every $(S,V)$ as in (1) of Theorem \ref{main},
there exists a distinguished triangulation $T$ such that every $(T,b)$
and $(T,b')$ can be explicitly connected by a sequence of inversions
of untrapped ambiguous edges.

(ii) For every $(S,V)$ as in (2) of Theorem \ref{main}, there exists a
distinguished triangulation $T$ such that for every $(T,b)$ and
$(T,b')$ either $(T,b)$ and $(T,b')$ or $(T,-b)$ and $(T,b')$ can be
explicitly connected by a sequence of inversions of untrapped
ambiguous edges.
\end{teo}

For every $(S,V)$ let us said {\it inversive} any triangulation $T$
which verifies the conclusions of Theorem \ref{A-main}, in accordance
to the two cases. Finally we have: 

\begin{teo}\label{inversive} For every $(S,V)$, every triangulation $T$
without trapped edges is inversive.
\end{teo}
\smallskip

{\it Proof of Theorem \ref{main-A}.} 
For every $S$ there is a minimum $n_S$ such that $\chi(S)- n_S<0$.
The proof is by induction on $n\geq n_S$. 
\medskip

{\bf Initial step: $(S,n_S)$.}

$\bullet$ $\ (\PP^2(\R),2)$. We use the naked triangulation $T$ of
Figure \ref{b-esiste}. Let $(T,b)$ and $(T,b')$ be supported by $T$.
By total inversion we can assume that $b$ and $b'$ agree on the
internal edges so that the internal vertex is a pit. The boundary
edges of the nutshell lift an ambiguous edge of both $(T,b)$ and
$(T,b')$; if $b\neq b'$ we conclude by inverting it in $b$. Then
$\tilde \Bb^{id}(\PP^2(\R),V)$ consists of one point.

\begin{prop}\label{PP-2} $|\Bb^{id}(\PP^2(\R),2)|=2$.
\end{prop}
\Dim Let $(T,b)$ and $(T,b')$ be as above such that $\delta(b,b')$
consists of the two internal edges. If we flip an internal edge of
$(T,b)$ we produce a trapped edge; in order to get the same naked
configuration we must flip the same edge in $(T,b')$ and $|\delta|$ is
unchanged. If we flip the edge of $(T,b)$ which lifts to the boundary
edges of the nutshell, we get a triangulation $(T_1, b_1)$ which is
abstractly like $(T,b')$ with respect to {\it another} nutshell, but
the two vertices exchange their role, so $(T_1,b_1)$ cannot be
relatively isotopic to $(T,b')$.

\cvd

\medskip

$\bullet$ $(S^2,3)$. Take $T$ made by two triangles glued along the
common boundary as in the first row of Figure \ref{b-esiste}.  Every
branched $(T,b)$ is determined by a labelling of the vertices by
$0,1,2$. Fix a $(T,b_0)$, then all $(T,b)$ are indexed by the elements
$\sigma$ of the symmetric group $\Sigma_3$, so that $(T,b_0)$
corresponds to the identity.  This group is generated by the
transpositions $\{(01),(12)\}$.  If $b=b_\sigma$, write $\sigma$ as a
product of minimal number of these generators.  Every such a sequence
of transpositions corresponds to a sequences of inversions of
ambiguous edges going from $(T,b_0)$ to $(T,b)$.

\medskip

In all other cases $n_S=1$.

$\bullet$ Let $S= S^1\times S^1$.
Consider $T$ as in the second row of Figure \ref{b-esiste}.
Every branching of $T$ is uniquely encoded by a total order
of the three vertices of one of the abstract triangles of $T$.
Then we can manage similarly as for $(S^2,3)$
by checking that also in this case every transposition
in a product of the generators corresponds to the inversion of 
an ambiguous edge.

\medskip

$\bullet$ Let $S$ be a Klein bottle. Refer to Figure \ref{b-esiste}.
If we use the triangulation $T$ made by two truncated bigons, it is
immediate that it carries exactly two branchings say $(T,b)$ and
$(T,-b)$. If we use as $T$ the other triangulation, we see that it
carries four branchings, distributed into two pairs $\{(T,b),
(T,b')\}$, $\{(T,-b), (T,-b')\}$ such that $(T,b')$ is obtained from
$(T,b)$ via the inversions of an ambiguous edge.

\medskip

Let us face now the remaining {\it generic cases} such that $\chi(S)<0$.

\medskip

{\bf Case (1)} Let $S$ be either orientable with $\chi(S)<0$ or non
orientable with $\chi(S)<0$ and even. Take the naked triangulation $T$
depicted in the proof of Lemma \ref{b-exist}; let $(T,b_0)$ be a
branched triangulation contructed therein. Let $(T,b)$ any branched
triangulation supported by $T$.  This determine a system of
orientations on the family of connection edges (Definition
\ref{connection-edge}).  Every connection edge is ambiguous in
$(T,b_0)$ so that up to some inversions, we can assume that $(T,b_0)$
and $(T,b)$ share such a partial system of orientations.  Let us cut
now $S$ along the connection edges. By restriction we get a family of
pairs of branched triangulated bricks $(B,b_0)$ and $(B,b)$ which
coincide at every connection edge. It is enough to show that every
$(B,b)$ is connected to $(B,b_0)$ by a sequence of inversions of
ambiguous edges. This can be checked case by case. We have three
types of $B$, the one-pierced torus or Klein bottle and the
twice-pierced torus. A priori, for every one-pierced brick we have two
local configurations of $(B,b_0)$ at the connection edges; for the
twice-pierced torus there are four. For every brick and every pair of
opposite local configurations, we see by means of the total inversion
that the desired result holds for one if and only if it holds for the
other. Then we are actually reduced to study one configuration in the
one-pierced cases, two in the twice-pierced one.  We organize the
discussion as follows, referring to Figure \ref{b-esiste2}:

\begin{itemize}
\item Denote by $t$ the top (abstract) triangle of $(B,b_0)$. Encode
 $(t,b_0)$ by labelling its vertices by $0,1,2$, say $v_0, v_1, v_2$;
 do the same for the bottom triangle $(t',b_0)$, getting $v'_0,v'_1,v'_2$.
 In the one-pierced brick $v_0=v'_0$.
 Every (abstract) edge of $(B,b_0)$ has two vertices belonging to 
 $\{v_0,v_1,v_2,v'_0,v'_1,v'_2\}$. For every $(B,b)$, every oriented
 edge $(e,b)$ will be denoted by its vertices, $e=ab$,
 written in the order so that the orientation emanates from the initial vertex $a$ 
 toward the final vertex $b$. For every one-pierced brick, we stipulate that the 
 connection edge in $(B,b_0)$ has $v_2$ as initial vertex.
 For the twice-pierced torus we stipulate that in $(B,b_0)$ 
 the pairs of connection edges is either 
 $(v_2v'_2, v_0v'_0)$ or $(v_2v'_2, v'_0v_0)$.
 We note that having fixed  orientation of the connection edges,
 $b_0$ is completely determined by $(t,b_0)$, that is this
 propagates in a unique way to a global branching.
  
  \item For every permutation $\sigma \in \Sigma_3$ we consider the 
  corresponding branched triagulated triangle $(t,b_\sigma)$ and we list 
  all the extensions to a global branching say $(B,b_\sigma)$, if any.
  Of course $(B,b_0)$ corresponds to the identity.
  
\item By varying $\sigma \in \Sigma_3$, the so obtained  $(B_\sigma,b)$ 
cover all possible branchings $(B,b)$ and we have
  to manage in order to connect $(B,b_0)$ with every $(B,b_\sigma)$.  It is
  convenient to start with the generating transpositions $\sigma =
  (0,1), (1,2)$, and express all other $\sigma$ as a product of three
  or two generators.
  \end{itemize}

Let us pass now to the actual verifications.
\smallskip

{\bf The one-pierced torus.}

$\sigma= (0,1)$: there a unique extension $(B,b_{(0,1)})$ which differs from
$(B,b_0)$ by the inversion of the ambiguous edge $v_1v_0$.
\smallskip

$\sigma = (1,2)$: there are several extensions. There is only one
containing the edge $v_0v_2$ and this differs from $(B,b_0)$ by the
inversion of the ambiguous edge $v_2v_1$. There are two extensions
containing the edge $v_2v_0$ which differ from each other by the
inversion of the ambiguous edge $v_0v'_2$. In the one containing
$v_0v'_2$, $v_2v_0$ is ambiguous, hence by inverting it we are in the
first case.

\smallskip

$\sigma = (0,1,2)=((1,2)(0,1)$: there are two extensions which differ
from each other by the inversion of the ambiguous edge $v_0v'_2$. In
the one containing $v_0v'_2$, the edge $v_0v_2$ is ambiguous, hence
possibly by inverting it we reach the case $(B,b_{(1,2)})$.

\smallskip

$\sigma=(0,2,1)=(0,1)(1,2)$: the discussion is similar to the one for
$(0,1,2)$; up to some inversion of ambiguous edges we reach the case
$(B,b_{(0,1)})$.

\smallskip

$\sigma = (0,2)=(0,1)(1,2)(1,0)$: there only one extension in which
$v_0v'_1$ is ambiguous.  By inverting it we reach the case
$(B,b_{(0,2,1)})$.

The first verification is complete.

\smallskip

{\bf The one-pierced Klein bottle.}

$\sigma= (0,1)$: there a unique extension $(B,b_{(0,1)})$ which
differs from $(B,b_0)$ by the inversion of the ambiguous edge
$v_0v_1$.
\smallskip

$\sigma = (1,2)$: there are no extensions.

\smallskip

$\sigma = (0,1,2)$: there is only one extension $(B,b_{(0,1,2)})$. The
following sequence of inversions of ambiguous edges realizes a transit
from this extension to $(B,b_{(0,1)})$ (we indicate the initial
orientation before the inversion): $v_0v_1$, $v'_2v_0$, $v_2v_1$,
$v_2v_0$.

 \smallskip

$\sigma=(0,2,1)$: there is only one extension. The following sequence
 of inversions of ambiguous edges realizes a transit to $(B,b_0)$:
 $v_2v_0$, $v_1v_0$.

\smallskip

$\sigma = (0,2)$: there are two extensions which differ to each other
by the ambiguous edge $v_0v'_2$.  In the ones containing the oriented
$v'_2v_0$, the edge $v_1v_0$ is ambiguous.  By inverting it we reach
$(B,b_{(0,1,2)})$.

\smallskip

The second verification is complete.
\smallskip

{\bf The twice-pierced torus.}
We have two cases depending on the orientation 
either $(v_2v'_2, v_0v'_0)$ or $(v_2v'_2, v'_0v_0)$
of the two connection edges.
\smallskip

{\bf Subcase  $(v_2v'_2, v_0v'_0)$}

$\sigma= (0,1)$: There are two extensions which differ by the
ambiguous edge $v'_0v'_2$.  In the ones containing the oriented edge
$v'_0v'_2$, $v_1v_0$ is ambiguous and possibly inverting it we reach
$(B,b_0)$.

\smallskip

$\sigma = (1,2)$: this is very similar to the case $(0,1)$.   

\smallskip

$\sigma = (0,1,2)$: there are several extensions. There is only one
containing $v'_2v_0$ which is ambiguous. In the ones containing
$v_0v'_2$ both $v_0v_2$ and $v'_0v'_2$ are ambiguous, then after at
most two inversions we reach $(B,b_{(1,2)})$.

\smallskip  

$\sigma=(0,2,1)$: this is very similar to the case $(0,1,2)$; via a
sequence of inversions we reach now $(B, b_{(0,1)})$.

\smallskip

$\sigma = (0,2)$: there are four extensions which differ by suitable
inversions of the edges $v_1v_2$ and $v_0v'_2$ which are both
ambiguous. Then up to such inversion we reach $(B,b_{(0,1)})$.

\smallskip

{\bf Subcase $(v_2v'_2, v'_0v_0)$}

At this point the fourth verification 
is a routine, we leave it to the reader.

 \medskip
 
 {\bf Case (2)} Let $S$ be not orientable such that $\chi(S)<0$ and
 odd.  We manage as in Case (1). The only difference is that the
 capping pierced Klein bottle is replaced with a pierced projective
 plane.  Again we use the triangulations depicted in the proof of
 Lemma \ref{b-exist}.  Up to total inversion we can assume that
 $(T,b_0)$ and $(T,b)$ coincide on the capping truncated bigon. The
 rest of the proof is unchanged.

The proof of Theorem \ref{main-A}  for $(S,n_S)$ is now complete.
\medskip

{\bf The inductive step.} Let us face first the generic case
$\chi(S)<0$.  So we have proved the result for $(S,1)$, and we want to
prove it for every $(S,n)$ by induction on $n\geq 1$.  We define the
distinguished triangulation for $(S,n)$ by modifying the one used for
$(S,1)$ as follows:

\begin{itemize}
\item The pierced Klein bottle and the twice-pierced torus bricks are unchanched.
\item We modify only the the one-pierced torus brick, say $B_1$, used
  when $n=1$ in order that $B_n$ carries all further $n-1$
  vertices. We do it inductively as follows: $B_1$ is triangulated by
  say $T_1$ as above; the naked triangulation $T_n$ of $B_n$ is
  obtained from $T_{n-1}$ by performing a $1\to 3$ move on the
  triangle which contains the connection edge.
\item In the treatment of $n=1$ we have also indicated a reference
  branched brick $(B,b_0)$; we define inductively the reference
  branching $(B_n,b_n)$ for every $n\geq 1$ as follows: set
  $(B_1,b_1):=(B,b_0)$ and recall that it is completely determined by
  a suitable total order of the vertices of the top triangle (labelled
  by $0,1,2$); $(B_n,b_n)$ is uniquely determined by extending the
  ordered set of vertices which defines $b_{n-1}$ by adding the new
  vertex produced by the $1\to 3$ move and stipulating that it is the
  smallest one (the vertices are labelled by $0,1,2, \dots, n-1, n,
  n+1, n+2$ and the labels of the vertices relative to $b_{n-1}$ shift
  by one).
\end{itemize}

\smallskip

We can fix the orientation of the connection edge, say
$v_{n+2}v'_{n+2}$. Consider any $(B_n,b)$.  There are two
possibilities:

\smallskip

(a) The new vertex of $B_n$ with respect to $B_{n-1}$ has a good
triangular star in $(B_n,b)$.  Then $b$ restricts to a branching
$(B_{n-1},b)$. By induction, this is connected to $(B_{n-1},b_{n-1})$
by a sequence of inversions of ambiguous edges. Finally we conclude by
applying Lemma \ref {2t-star} to the innermost triangular star.
\smallskip

(b) The innermost triangular star as above is bad in
$(B_n,b)$. Consider first $(B_2,b)$; we readily see that $v_{2}$ is
necessarily either a pit or a source. In any case $v_nv_{n+2}$ is
ambiguous; by inverting it the triangular star becomes good and we
reach the case (a). In general we can assume by induction that the
triangular star of $v_1$ with respect to the restriction of $b$ to
$B_{n-1}$ is good, so that we can apply the above raisoning to the
triangular star of $v_0$ in $(B_n,b)$ and reach again the case (a).
\smallskip

The proof of Theorem \ref{main-A} in the generic cases is now complete.

\cvd

\smallskip

For the remaining cases such that $\chi(S) \geq 0$ we limit ourselves
to some indications.
\smallskip

$\bullet$ $(S^2,n)$, $n\geq 3$. Denote by $T_3$ the triangulation used
above for $n=3$. Select one triangle $t$ and one edge $e$. For every
$n>3$, the distinguished triangulation $T_n$ for $(S^2,n)$ is obtained
by induction on $n$ by performing a $1\to 3$ move on the triangle of
$T_{n-1}$ which is contained in $t$ and contains $e$.  In particular
$T_4$ corresponds to the triangulation of the boundary of a
tetrahedron. Every $(T_4,b)$ is determined by a labelling of the
vertices by $0,1,2,3$.  Fix a $(T,b)$; then the branchings are indexed
by the elements of the symmetric group $\Sigma_4$. This is generated
by the transpositions $(0,1),(1,2),(2,3)$.  Write every $\sigma$ as a
product of these generators with minimal number of terms. This
corresponds to a sequence of inversions of ambiguous edges connecting
$(T,b)$ and $(T,b_\sigma)$. For $n>4$ we argue by induction on $n$.
\smallskip

$\bullet$ $(S^1\times S^1,n)$ or $(\PP^2(\R) \# \PP^2(\R),n)$, $n\geq
1$. In both cases we start with the triangulation say $T_1$ used for
$n=1$. Then (referring to Figure \ref{b-esiste}) $T_n$ is obtained
from $T_{n-1}$ by performing a $1\to 3$ move on the triangle contained
in the top triangle and containing the diagonal edge of $T_1$.
\smallskip

$\bullet$ $(\PP^2(\R),n)$, $n\geq 2$. We start with $T_2$ used for
$n=2$. Referring to Figure \ref{b-esiste}, $T_3$ is obtained by
performing a bubble move at the internal edge on the left side. Denote
by $t$ the new triangle contained in the top half of $T_2$ and by
$e$ its edge contained in the interior of this top-half. Then, for
$n>3$, $T_n$ is obtained from $T_{n-1}$ by performing a $1\to 3$ move
on the triangle of $T_{n-1}$ wich contains $e$.
\smallskip

The proof of Theorem \ref{main-A} is now complete.

\cvd

\smallskip

{\it Proof of Theorem \ref{inversive}.} 

{\bf Case (a):} $\chi(S)$ is not strictly negative and odd. The
distinguished inversive triangulations $T_n$ of $(S,n)$ constructed in
the proof of Theorem \ref{main-A} have no trapped edges. Let $T$ be
any other triangulation without trapped edges. We know that there is a
sequence of flips $T_n \Rightarrow T$ through triangulations without
trapped edges. Denote by $l$ the number of flips. We work by induction
on $l$. Let $T'$ be obtained by performing the first $l-1$ flips. By
induction the theorem holds for $T'$. Hence we are reduced to check
the case $l=1$. We fix the notations as follows: $e$ is the flipping
edge, that is a diagonal of a quadrilateral $Q=t_1\cup t_2$ in $T_n$,
$t_1\cap t_2=e$; $e'$ denote the other diagonal of $Q$, that is the
edge of $T$ which replaces $e$.  Let $(T_n,b)$ and $(T_n,b')$ be
connected by a sequence of $k$ inversions of ambiguous egdes. Let
$(T,\tilde b)$, $(T,\tilde b')$ be obtained by $b$-enhancing in some
way the flip $(T_n,b)\to T$ and $(T_n,b')\to T$. We want to modify the
sequence in order to get one connecting this branchings of $T$. We
note that if a $b$-flip is not forced then the two possibilities are
related by inverting an ambiguous edge, so this is essentially
immaterial for our discussion. If an inversion concerns an edge not
contained in $(Q,b)$ then it makes sense also on $(T,\tilde b)$. By
these remarks and working by induction on $k$, we are reduced to
analyze the inversion of an edge $e^*$ contained in $(Q,b)$. There are
a few possibilities.
\smallskip

$\bullet$ $e^*=e$; then the $b$-flip is either totally ambiguous
(i.e. bump) or forced ambiguous, depending if the vertices of $(Q,b)$ opposite
to $e$ are either both a pit (resp. source) or one is a pit and the
other a source.  So in the first case we possibly replace the
inversion of $e$ with the inversion of $e'$.
\smallskip

Assume now that $e^*\neq e$.

$\bullet$ $e$ is ambiguous in $(Q,b)$ as above. If the flip is bump,
there are two possibilities for $e^*$.  Then the inversion of $e^*$
can be performed on $(T,\tilde b)$, possibly after having inverted the
ambiguous edge $e'$. If the flip is forced ambiguous, again there are
two possibilities for $e^*$ and in every case the inversion of $e^*$
can be performed on $(T,\tilde b)$.

\smallskip

$\bullet$ $e$ is ambiguous in one of the two triangles of $(Q,b)$ and
non ambiguous in the other. The flip is non ambiguous. There are three
possibilities for $e^*$. We readily check that in every case the
inversion of $e^*$ can be performed on $(T,\tilde b)$.

\smallskip

$\bullet$ $e$ is non ambiguous in both triangles of $(Q,b)$. The flip
is not forced. We check that the inversion of $e^*$ can be performed
on $(T,\tilde b)$, possibly after having inverted the ambiguous edge
$e'$.
\smallskip

This complete the proof in Case (a).

\smallskip

{\bf Case (b):} $\chi(S)$ is strictly negative and odd. The
distinguished triangulations $T_n$ of $(S,n)$ have one trapped edge
carried by the one-pieced projective plane.  Let $T^*_n$ be obtained
by flipping its connection edge. $T^*_n$ does not contain any trapped
edge and arguing similarly as above, we see that it is inversive. Then
the proof is like in Case (a), by using $T^*_n$.
\smallskip

Theorem \ref{inversive} is achieved.

\cvd     

\section{B-proof of  the main Theorem} \label{B-main}  
We are going to implement strategy {\bf B}. For simplicity we will deal only
with the generic case $\chi(S)<0$. Strictly speaking
Remark \ref{AvsB} will apply to the case of orientable $S$,
as in the non orientable case we will actually adopt a mixture
of {\bf A} and {\bf B}.

 \subsection{B-proof when $S$ is orientable} \label{B-or}  
Let $(T,b)$ and $(T,b')$ triangulations of $(S,V)$ such
that $\delta(b,b')$ is non empty.  By using Lemma \ref{no-trapped} and
Lemma \ref{eliminate-amb} it is not restrictive to deal under the following:

\medskip

{\bf Initial assumptions:} {\it (1) $T$ does not contain trapped
  edges;

(2) Every disoriented edge $e\in \delta(b,b')$ is {\it non} ambiguous
  in both $(T,b)$ and $(T,b')$.}

\medskip  

At first we analyze the effects of {\it flipping $e\in \delta(b,b')$}
in both $(T,b)$ and $(T,b')$ looking for a decreasing of $|\delta|$ if
any. Let $t_1$ and $t_2$ be the two triangles of $T$ which share $e$.
As $e$ is non ambiguous in both triangulations, then $e$ is non
ambiguous in at least one of the branched triangles $(t_j,b)$ and
similarly for the $(t_j,b')$'s. There are two possibilities:
\begin{enumerate}

\item There is at least one triangle, say $t_1$, such that $e$ is non
  ambiguous in both $(t_1,b)$ and $(t_1,b')$, so that necessarily
  $(t_1,b')=(t_1,-b)$.

\item $e$ is non ambiguous (resp. ambiguous) in $(t_1,b)$
  (resp. $(t_2,b)$) while $e$ is non ambiguous (resp. ambiguous) in
  $(t_2,b')$ (resp. $(t_1,b')$).
\end{enumerate}

\medskip

Let us analize the first case.
\medskip

{\bf Case (1).} Concerning $t_2$, there are three possibilities: it
can contains either $k=0,1$ or $2$ further edges belonging to
$\delta(b,b')$. Note that $k=2$ if and only if $b'=-b$ on the whole of
$t_1\cup t_2$. We say that the disoriented $e$ is {\it (1)bad} if
$k=2$ and $e$ is ambiguous in $(t_2,b)$ (hence in $(t_2,b')$). In all
other cases we say that $e$ is {\it (1)good}. We have

\begin{lem}\label{1goodOK} 
  Let $e\in \delta(b,b')$ be (1)good. Then 
  by flipping $e$ in both $(T,b)$ and $(T,b')$
  we get $(T',b_1)$ and $(T',b_2)$ which still satisfy our assumptios
  and such that $|\delta(b_1,b_2)|<|\delta(b,b')|$.
\end{lem}
\Dim  First we note that if $f_e$ creates a trapped edge, then
necessarily $e$ is internal to some naked nutshell say $N$ in $T$.
Now  we analize the situation case by case according to the value of
$k=0,1,2$.

\begin{figure}[ht]
\begin{center}
 \includegraphics[width=8cm]{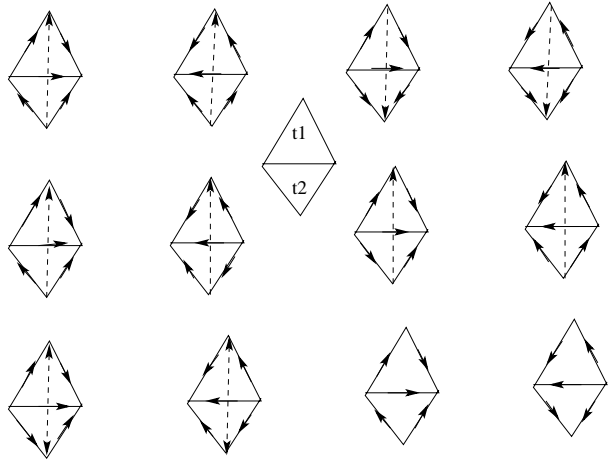}
\caption{\label{untrapped1} Flipping an untrapped edge, Case (1).}
\end{center}
\end{figure}

\begin{figure}[ht]
\begin{center}
 \includegraphics[width=7cm]{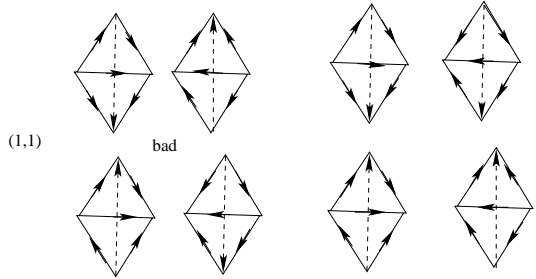}
\caption{\label{untrapped2} Flipping an untrapped edge, Case (2).}
\end{center}
\end{figure}

\begin{enumerate}
\item If $k=0$  both $f_{e,b,b_1}$ and $f_{e,b',b_2}$ are
non ambiguos and $|\delta|$ decreases by $1$. There are
not compatibly branched nuthshells containing $e$, hence
the flip does not create any trapped edge (see the first row
  of Figure \ref{untrapped1}).

\item If $k=1$,  we can choose $f_{e,b,b_1}$ and $f_{e,b',b_2}$ 
in such a way that $|\delta|$ decreases by $1$ and the
 new edge is ambiguous in one of the triangulations obtained so far
 (see Figure \ref{untrapped1}, second row). We claim that $f_e$
 does not  create any trapped edge. Otherwise one among 
 the branched nutshells $(N,b)$ and $(N,b')$ would be good with
 internal vertex which is either a pit or a source. Then (recall  
 Lemma \ref{2_nutshell} (3)) it would contain an  ambiguous internal edge
 $\tilde e$ belonging to $\delta(b,b')$ against our assumption.
 
\item If $k=2$ and $e$ is (1)good, then necessarily $e$ is non
  ambiguous also in $(t_2,b)$ (hence in $(t_2,b')$).  
  We can choose $f_{e,b,b_1}$ and $f_{e,b',b_2}$  in such a way that
  $|\delta|$ decreases by $1$ and the new edge is ambiguous in
  both triangulations obtained so far (see the left side of the third row
  of Figure  \ref{untrapped1}).  If $e$ would belong to a
  naked nutshell in $T$, then  both $(N,b)$ and $(N,b')$ are good
  with internal vertex which is either a pit or a source, and we can
  argue as above. So $f_e$ does not create any trapped edge.

\end{enumerate}

The Lemma is proved.

\cvd
 
Concerning the $(1)bad$ situation, we readily realize that  

\begin{lem}\label{1bad} 
  If $e$ is (1)bad, then by flipping $e$ we keep the same value of
  $|\delta|$.  The new edge is non ambiguous in both triangulations
  obtained so far (see the right side of the third row of Figure
  \ref{untrapped1}).  If the flip creates a trapped edge, then both
  nuthshells $(N,b)$ and $(N,b')$ are bad and all edges of $N$ belong
  to $\delta(b,b')$
\end{lem}   

\medskip

Let us turn now to the second case.

\medskip

{\bf Case (2).}   
It is easy to see that necessarily 
the boundary of $t_1\cup t_2$
contains exactly a couple $e_1\subset t_1, \ e_2\subset t_2$ of edges  
which do not belong to $\delta(b,b')$.
There are two possibilities, see Figure \ref{untrapped2}.  

\medskip

(i) $e_1$ and $e_2$ are consecutive edges in the (abstract)
quadrilateral $t_1\cup t_2$.  In this case we can flip $e$ and
$|\delta|$ decreases by 1. If this creates a trapped edge, then as
usual $e$ is an internal edge of a nutshell $N$, and both $(N,b)$ and
$(N,b')$ are bad with the two boundary edges which do not belong to
$\delta(b,b')$.

\medskip

(ii) $e_1$ and $e_2$ are opposite edges in the (``abstract")
quadrilateral $t_1\cup t_2$ and their orientations are necessarily compatible,
that is they extend to an orientation of the whole boundary of the
quatrilateral.  Then by flipping $e$ we keep the same value of
$|\delta|$. The flip does not create any trapped edge.

\medskip

If we are in case (i) and the flip does not create a trapped edge,
then we say that $e$ is {\it (2)good}. The other cases are {\it (2)bad}.
Let us call generically {\it good} an edge $e\in \delta(b,b')$ which
is either (1)good or (2)good.  Otherwise let us say that it is {\it bad}.
Summarizing the above discussion we can strenghten the ``Initial assumptions''
as follows:

\begin{lem}\label{all-bad} In order to prove Theorem
\ref{main} it is not restrictive to deal under the following
\smallskip

{\bf All-bad assumptions:} (1) $T$ does not contain trapped
edges; 

(2) Every edge $e\in \delta(b,b')$ is {\it non} ambiguous in both $(T,b)$ and
$(T,b')$;

(3) Every edge in $\delta(b,b')$ is bad.
\end{lem}

 \medskip

 Now we show that the all-bad
 assumptions are quite constraining.
 
 \begin{lem}\label{all-2bad} Let $(T,b)$ and $(T,b')$ satisfy
the all-bad assumptions. Then every $e\in \delta(b,b')$
is actually (2)bad.
\end{lem}
\Dim Assume that there exists a (1)bad edge $e\in \delta(b,b')$.  This
propagates to all edges of $T$: eventually $b'=-b$ and all edges should be
(1)bad. We want to show that this is impossible.  We note that there
are not adjiacent (necessarily bad) nuthshells because a common
boundary edge should be ambiguous, hence good.  
Let $v$ be a vertex of
$T$ which is not the center of a nutshell and analyse the
possible configurations of its (abstract)  developed {\it star} $St(v,b)$ 
$(T,b)$ (the one in $(\tilde T,b')$ is obtained by just reversing the
orientations).

{\bf Claim:} Every such a star $St(v,b)$ has the following qualitative
configuration. Every edge in the boundary of the star is {\it
  ambiguous} in the respective triangle. $St(v,b)$ can contains an
{\it even} number of bad nutshells sharing the vertex $v$ 
(necessarily even because otherwise the star would contain a good edge) 
and the orientations of their
boundaries alternate (with respect to the reference orientation of
$\tilde S$). The edges of the boundary of the star between two consecutive nutshells
have compatible orientations as well as the internal edges of their respective 
triangles are all either
ingoing or outgoing with respect to the central vertex $v$.  Moving
along the boundary of the star, boundary orientations and ``in-out''
types switch each time we pass a nutshell. In particular if there are
no nutshells, then the boundary of $St(v,b)$ is an oriented circles.
\medskip

{\it Proof of the Claim.}  Assume that there is an edge $e$ in the
boundary of $St(v,b)$ which is {\it non} ambiguous in the relative
triangles. Let us try to complete the star by moving along its
boundary in the direction of the orientation of $e$. Possibly after
some boundary edges which are oriented like $e$ and are {\it
  ambiguous} in the respective triangle (with internal edges pointing
towards the central vertex $v$), we necessarily find either a boundary
edge $e'$ which is non ambiguous in the respective triangle and has
opposite orientation with respect to $e$, or a bad nutshell (whose
boundary orientation is uniquely determined). We see that in both case
there is an internal edge which is ambiguous in the star, hence good.

\cvd

Now we can conclude by noticing that for every $St(v,b)$ with the
properties stated in the Claim there is a vertex $v'$ in the boundary
of the star such that the boundary of $St(v',b)$ contains an edge
which is non ambiguous in the relative triangle of $St(v',b)$.  
Lemma \ref{all-2bad} is proved.

\cvd

\begin{remark}\label{orientable}
  {\rm The hypothesis that $S$ is orientable has been already employed
    in order to limit the way a trapped edge can be produced; it will
    be important also in the rest of the discussion; the key point is
    that it prevents that the stars of the disoriented (2)bad edges
    (in a all-bad configuration) glue each other at edges not
    belonging to $\delta(b,b')$ producing a M\"obius strip. For
    example, the opposite edges not belonging to $\delta(b,b')$ in a
    basic (2)bad configuration $(t_1\cup t_2,b) $, $(t_1\cup t_2,b')$
    cannot be identified.  }
\end{remark}
     
\begin{defi}\label{terminal} {\rm 
(1) A {\it terminal (2)bad type} is either:
\begin{itemize}
\item   A couple of (2)bad nutshell $(N,b)$, $(N,b')$ such that  
   $|\delta(b,b')|=2$ and the boundary edges do not belong to $\delta(b,b')$.
   \item A couple of triangulated annuli $(A,b)$, $(A,b')$ obtained from
    a basic (2)bad  configuration $(t_1\cup t_2,b) $, $(t_1\cup
    t_2,b')$ by identifying the opposite boundary edges of the
    quadrilateral which belong to $\delta(b,b')$. For the resulting
    triangulations we have $|\delta(b,b')|=2$ and the boundary
    of $(A,b)$ (similarly for $(A,b')$) is formed by two circles each
    containing one vertex and endowed with opposite orientations. 
    \end{itemize}
    
    (2) A couple $(T,b)$ $(T,b')$ is said {\it terminal all-bad} if
verify the all-bad assumptions and every disoriented edge $e\in
\delta(b,b')$ is contained in a terminal (2)bad type.  }
\end{defi}

\begin{figure}[ht]
\begin{center}
 \includegraphics[width=3.5cm]{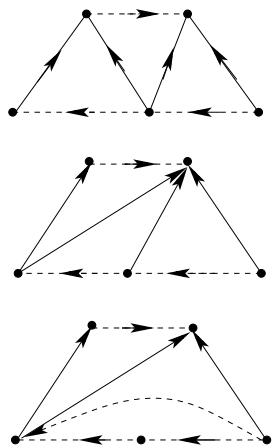}
\caption{\label{terminal} Terminal move.}
\end{center}
\end{figure}

We have

\begin{lem}\label{terminal-all-bad}
In order to prove Theorem
\ref{main} it is not restrictive to deal 
with {\bf terminal all-bad couples $(T,b)$, $(T,b')$}.
\end{lem}

\Dim Let $(T,b)$, $(T,b')$ verify the all-bad assumptions.  If $(T,b)$
presents a pattern as in the top of Figure \ref{terminal} (we
stipulate that the dashed edges do not belong to $\delta(b,b')$) we
can perform the sequence of $b$-flips suggested by descending the rows
of the picture (the corresponding flips on $(T,b')$ are
understood). We eventually get $(T',b_1)$ and $(T',b_2)$ such that
$|\delta(b_1,b_2)|=|\delta(b,b')|-1$ and the number of (2)bad $t_1\cup
t_2$ decreases by 1.  We stop when we reach a terminal configuration.

\cvd

We can state now the conclusive lemma.

\begin{lem}\label{final}
  Let $(T,b)$ and $(T,b')$ be a terminal all-bad couple.  Then we can
  find sequence of inversions of ambiguous edges $(T,b)\Rightarrow
  (T,\tilde b)$, $(T,b')\Rightarrow (T,\tilde b')$ such that
  $|\delta(\tilde b, \tilde b')|< |\delta(b,b')$.
\end{lem} 

By iterating all the above procedure starting from $(T,\tilde b)$,
$(T, \tilde b')$ we eventually get $|\delta|=0$ and the main theorem
follows.

{\it Proof of Lemma \ref{final}.}  Let $e$ be an edge not belonging to
$\delta(b,b')$ and contained in the star of a disoriented $\bar e \in
\delta(b,b')$. If $e$ is ambiguous, let us invert it in both $(T,b)$
and $(T,b')$.  If $e$ was in the boundary of a bad nutshell, after the
inversion the nuthshell becomes good and we can apply Lemma
\ref{2_nutshell}. If $e$ was in the boundary of a basic (2)bad
configuration $(t_1\cup t_2,b) $, $(t_1\cup t_2,b')$, then after the
inversion $\bar e$ becomes ambiguous (recall Remark \ref{orientable})
and we can invert it in $(T,b)$ to decrease $|\delta|$ by $1$.  So, if
$e$ is ambiguous we have done. Assume that $e$ is not ambiguous. Then
$e$ is the edge of a (abstract) triangle which is enterely formed by
edges not belonging to $\delta(b,b')$. Let $v$ be the vertex of this
triangle which does not belong to $e$.  We realize that there is an
internal edge of $St(v,b)$ which is ambiguous. Then by succesive
inversions of ambiguous edges we eventually make $e$ ambiguous an we
can conclude as above.

\cvd

\smallskip

The {\bf B}-proof in the orientable case is complete.

\subsection{AB-proof in the non orientable case}\label{AB-main}
In the non orientable case we do not use as initial situation an arbitrary
couple of triangulations (without trapped edges) $(T,b)$, $(T,b')$ of $(S,V)$.
We partially specialize the choice by requiring that $T$ respects a decomposition
$S= S_0  \#  N$ where $S_0$ is orientable and $N$ is either a Klein bottle
or a projective plane. More precisely we require that $T$ is the union 
of a one-pierced non orientable brick with the distinguished triangulation 
already used in the proof of Theorem \ref{main-A}, and an {\it arbitrary} 
triangulation of a one-pierced $S_0$, provided that they coincide at the connection
loop of the non orientable brick. Then the orientable $B$-proof applies with minor  
changes to the one-pierced $S_0$.

\cvd

\section{On the ideal sliding equivalence}\label{sliding}
The sliding transits have been already considered in Section 5 of \cite{NA}.
There we have been mainly concerned with the completed sliding
equivalence, more precisely with triangulations of a given oriented surface with 
arbitrary number of  vertices,  considered all together. 
Here we stress the ideal set up, we consider also non orientable surfaces
and we introduce some new constructions (for instance the so called 
horizontal foliation). The basic idea 
is that every branched triangulation of $(S,V)$
carries  some remarkable structures of differential topological type which 
are preserved by the ideal sliding  while they can be widely modified
by the bump transits accordingly with Theorem \ref{main}.

\subsection{The tranverse foliations carried by a branched triangulations} 
\label{carried} Given $(S,V)$ we recall that $S_V$ denotes the surface
with boundary obtained by removing a family of disjoint open $2$-disks
centred at each $v\in V$. We are going to consider possibly singular
foliations on $S_V$ or $S$. Every such a foliation can be obtained by
integration of some field of tangent direction (a tangent vector field
if the leaves are oriented). We will say that two foliations are {\it
  homotopic} ({\it isotopic}) if they are obtained by integration of
homotopic (isotopic) fields. Let $(T,b)$ be a branched triangulation
of $(S,V)$. First we are going to show that $(T,b)$ carries in a
canonical way a pair of regular transverse foliations $(\Vv,\Hh)$ on
$S_V$, called respectively the {\it vertical} and the {\it horizontal}
foliations. $\Vv$ is always oriented. If $S$ is oriented then also
$\Hh$ can be oriented in such a way that every intersection point has
intersection number equal to $1$. $\Vv$ and $\Hh$ can be extended to
singular foliations $\vG$ and $\hG$ defined on the whole of $S$. They
share the singular set $Z$ which consists of the vertices $v\in V$
where the index $1-d_b(v)\neq 0$.  The two foliations are transverse
on $S\setminus Z$.  If $S$ is oriented, then both $\vG$ and $\hG$ are
oriented; if the singularity indices are all non positive, so that $S$
is of genus $g\geq 1$, then $(\vG, \hG)$ looks like the couple of
vertical and horizontal foliations of the square of an Abelian
differential on a Riemann surface.  Let us pass to the actual
definition.

\begin{figure}[ht]
\begin{center}
 \includegraphics[width=7.5cm]{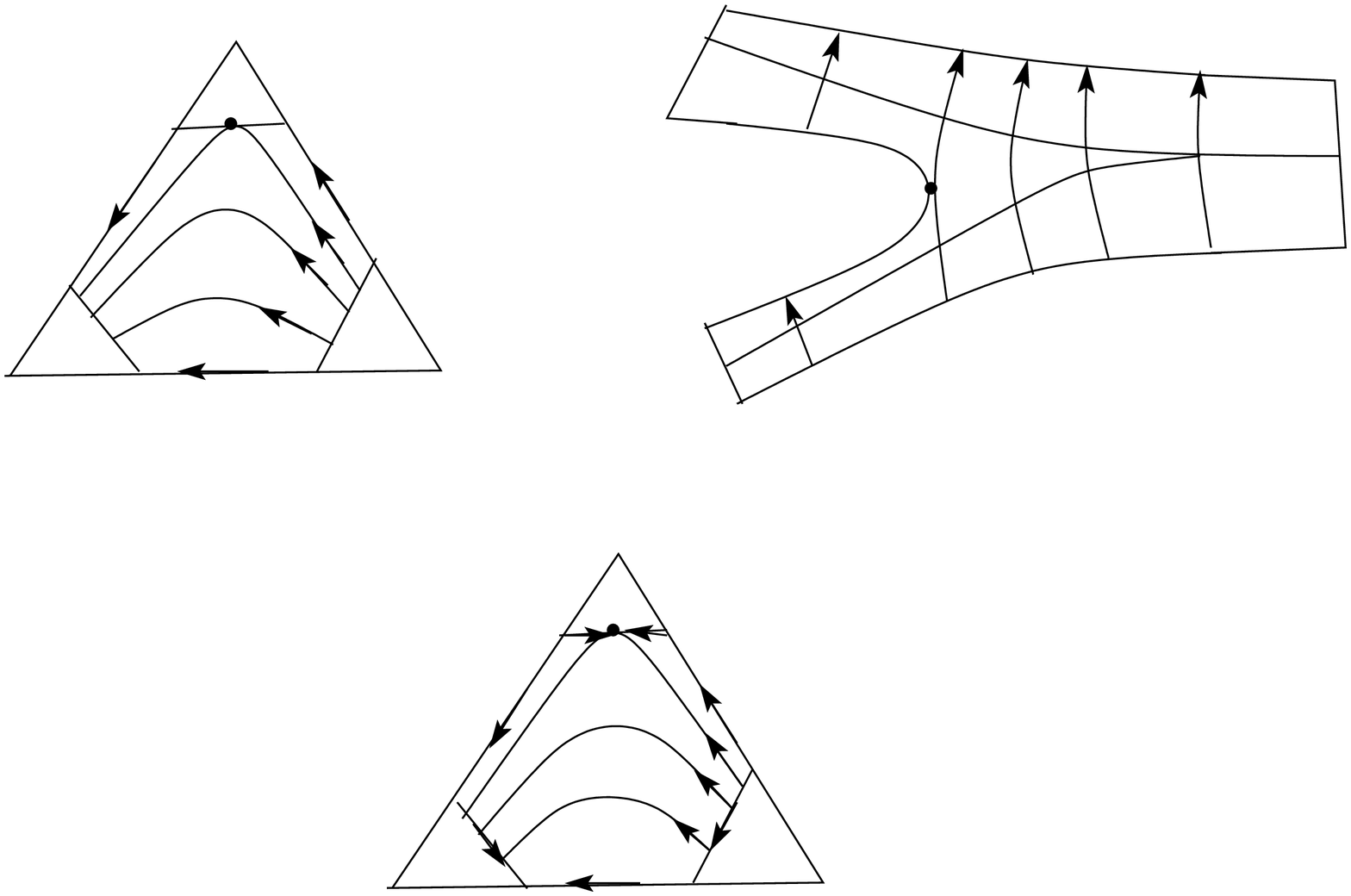}
\caption{\label{V} The vertical foliation on $S_V$.} 
\end{center}
\end{figure}

{\bf The vertical foliation $\Vv=\Vv(T,b)$.} Let $(T,b)$ be as above.
Up to isotopy, the intersection of $S_V$ with every (abstract)
triangle $t$ of $T$ is a ``truncated triangle'' $\bar t$, i.e. a
hexagon with $3$ internal ``long'' edges (each one contained in an
edge of $T$) and $3$ ``short'' edges contained in the boundary
$\partial S_V$.  The short edges are in bijection with the cornes of
triangles of $T$; some are labelled by $1$ in accordance with the
associated corners.  The union of the hexagons form a cell
decomposition of $S_V$, by the restriction to the long edges of the
gluing in pairs of the abstract edges of $T$.  The union of the short
edges form a triangulation of $\partial S_V$.  By using the branching
$b$, we are going to endow every hexagon $\bar t$ of an oriented
foliation $\Vv(\bar t,b)$. These constitute the tiles of a puzzle that
once assembled realizes $\Vv$. This is illustrated on the top-left of
Figure \ref{V}. In fact $\Vv(\bar t,b)$ is the restriction to $(\bar
t, b)$ of a classical Whitney field which can be defined explicitely
in terms of barycentric coordinates (see \cite{HT}).  In the dual wiew
point, recall that the spine $(\theta,b)$ of $S_V$ is an embedded
transversely oriented train track in $S_V$; the foliation $\Vv$ is
positively transverse to it (on the top-right of Figure \ref{V} we see
the dual picture corresponding to the puzzle tile).  The foliation
$\Vv$ has remarkable properties.
\begin{defi} {\rm A {\it traversing foliation} on
    $S_V$ is a foliation with oriented leaves such that:
\begin{enumerate}
\item Every leaf of $\Ff$ is a closed interval which
  intersects transversely $\partial S_V$ at its endpoints.
\item There is a {\it non empty} finite set of exceptional leaves of $
  \Ff$ which are simply tangent to $\partial S_V$ at finite number of
  points.
  \end{enumerate}
  
  $\Ff$ is {\it generic} if every exceptional leaf is tangent
  to the boundary at one point.
}
\end{defi}

Then it is not hard to see that:
\begin{prop}\label{traversing} 
(1) The vertical foliation $\Vv$ associated to a branched spine
  $(T,b)$ of $(S,V)$ is a generic traversing foliation on $S_V$. The
  exceptional leaves of $\Vv$ are in bijection with the $1$-labelled
  short edges of $\partial S_V$; every exceptional leaf is tangent to
  the interior of the associated edge.  $\Vv$ is uniquely determined
  up to isotopy.

(2) The dual branched spine $(\theta,b)$ of $S_V$ intersects
  transversely all leaves of $\Vv$.  Every exceptional leaf intersects
  transversely $\theta$ at two points. A leaf passing through a
  singular point of $\theta$ is generic and intersects $\theta$ at one
  point.  A generic leaf intersects $\theta$ at one or two points. In
  the second case it is contained in a quadrilatelar in $S_V$
  vertically bounded by an exceptional leaf and a leaf passing through
  a singular point of $\theta$.

(3) Every generic traversing foliation on $S_V$ can be realized as the vertical
foliation of some branched spine of $(S,V)$.
\end{prop}

{\bf Boundary bicoloring.} Given a traversing foliation $\Ff$ of
$S_V$, denote by $X=X_\Ff \subset \partial S_V$ the set of tangency
points of the exceptional leaves.  $\Ff$ determines a {\it bicoloring}
of the components of $\partial S_V \setminus X$, denoted by $\partial
\Ff$; let us say that a component $c$ is {\it white} (resp. {\it
  black}) if the foliation is ingoing (outgoing) along $c$. If $S$ is
oriented the color can be encoded by an orientation, in the sense that
a black component keeps the boundary orientation of $\partial S_V$
(according to the usual rule ``first the outgoing normal''), while a
white one has the opposite orientation.
In the bottom of Figure \ref{V} we see the oriented enhancement of the
$\Vv$-tile (we stipulate that in the picture the $b$-orientation of
the triangle agrees with the orientation of $S_V$; we obtain the
picture for the negative branching $-b$ by just inverting all arrows).

\smallskip

\begin{figure}[ht]
\begin{center}
 \includegraphics[width=7.5cm]{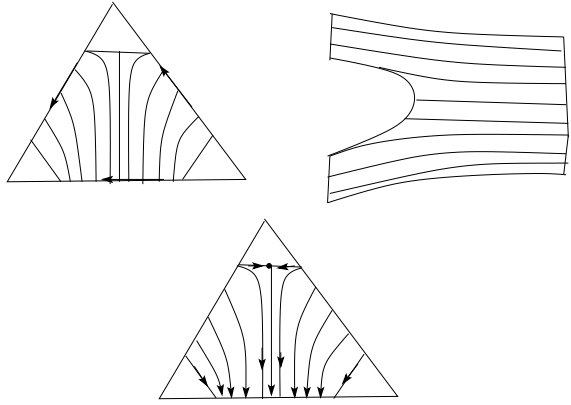}
\caption{\label{H} The horizontal foliation on $S_V$.} 
\end{center}
\end{figure}

{\bf The horizontal foliation $\Hh= \Hh(T,b)$.}  Alike $\Vv$ we define
$\Hh$ as the result of a puzzle. On the top of Figure \ref{H} we show
the foliated hexagon and the corresponding dual picture.  In general
$\Hh$ is not oriented. If $S$ is oriented then $\Hh$ is oriented as
well; on the bottom of Figure \ref{H} we show the oriented version of
the tile. Now we realize that:

\begin{enumerate}
\item $\Vv$ and $\Hh$ are transverse foliations.
\item Let $Y$ be the union of the $1$-labelled short edges.  Then
  every component of $\partial S_V \setminus Y$ is contained in a leaf
  of $\Hh$, while $\Hh$ is transverse to the interior of every
  $1$-labelled short edge.
\item If $S$, hence $\Hh$, is oriented then every boundary component
  in a leaf is oriented like the leaf and is contained in a component
  of $\partial S_V \setminus X$; these orientations propagate to the
  whole components of $\partial S_V \setminus X$ and reproduce the
  bicoloring orientation.  $\Vv$ intersects $\Hh$ with intersection
  number equal to $1$ everywhere.
\item The pair $(\Vv,\Hh)$ is uniquely defined up to isotopy.
\end{enumerate}

\begin{figure}[ht]
\begin{center}
 \includegraphics[width=6.5cm]{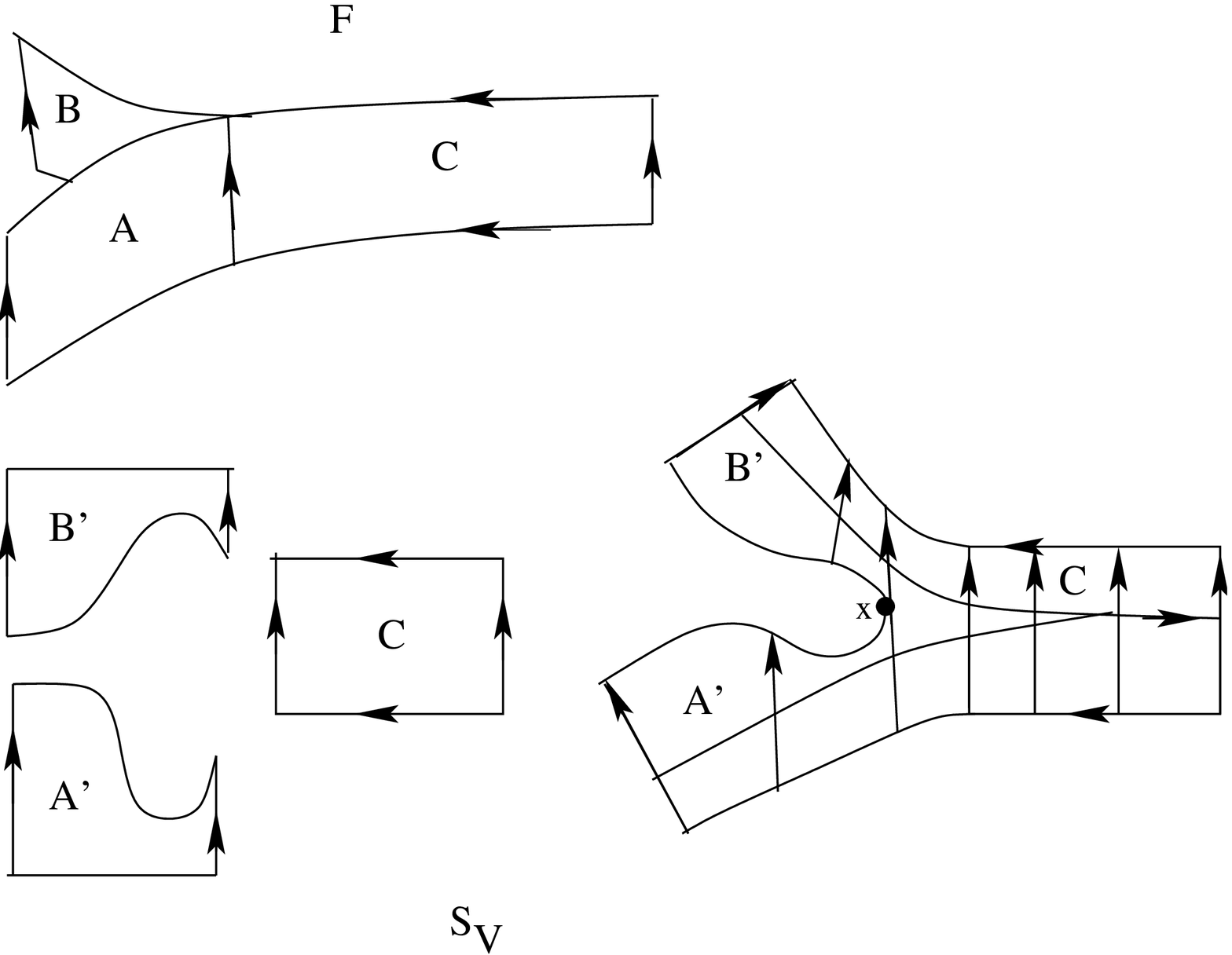}
\caption{\label{vertical} Another realization of $\Vv$.} 
\end{center}
\end{figure}

If $S$ is oriented there is another way to
    realize the foliations $(\Vv,\Hh)$ by using the 2D case of a
    result of \cite{GR}. This can be described as follows. Take a (non
    embedded) copy $\theta^*$ of the oriented train track
    $(\theta,b)$. Consider the oriented branched surface $F:= \theta^*\times
    [-1,1]$; this carries the verical foliation $\Vv^*$ with leaves of
    the form $\{x\} \times [-1,1]$ and the horizontal one $\Hh^*$ with
    branched leaves of the form $\theta^*\times \{y\}$.  Then one can
    find an embedding of $S_V$ into $F$ which preserves the orientation and 
    such that $\Vv$ is just the
    restriction of $\Vv^*$ to $S_V$. This is suggested in Figure
    \ref{vertical}. Also $\Hh^*$ restricts to a regular foliation of
    $S_V$ which becomes our final horizontal foliation $\Hh$ after a
    suitable homotopy.

\smallskip

{\bf Extension to singular foliations.} Let $(S,V)$ be as usual.
Let us define some notions.
A function
$$\iG: V\to \{n\in \Z | n\leq 1 \}$$
is {\it admissible} if  
$$\chi(S)= \sum_v \iG(v) \ . $$ 
For every admissible
function $\iG$, a {\it vertical foliation of type $\iG$} on $(S,V)$
is an oriented singular foliation $\vG$ that
verifies by definition the following properties:  

\begin{enumerate} 
\item The singular set $Z$ of $\vG$ consists of the $v\in V$ such that
  $i(v)\neq 0$.
\item If $\iG(v) \neq 0,1$, then the local model of $\vG$ at $v$ is
  given by the vertical foliation at $0$ of the quadratic differential
  $z^{-2i(v)}dz^2$. If $\iG(v)=1$, the local model is given by the
  integral lines at $0$ of the gradient of either the function $|z|^2$
  or $-|z|^2$.
\end{enumerate}

\smallskip

An {\it horizontal foliation of type $\iG$} on $(S,V)$
is a non oriented singular foliation $\hG$ that
verifies by definition the following properties:  

\begin{enumerate} 
\item  The singular set $Z$ of $\hG$ consists of the $v\in V$
such that $\iG(v)\neq 0$.
\item If $\iG(v) \neq 0,1$, then the local model of $\hG$ at $v$
is given by the horizontal foliation at $0$ of the quadratic differential
$z^{-2i(v)}dz^2$. If $\iG(v)=1$, the local model is given by the 
level curves at $0$ of the function $|z|^2$.
\end{enumerate}

A {\it transverse pair of foliations of type $\iG$} is a pair
$(\vG,\hG)$ such that 
\begin{enumerate}
\item $\vG$ and $\hG$ are vertical and horizontal foliations
of type $\iG$ respectively.
\item The two foliations are
transverse on  $S\setminus Z$ and, case by case, the above local models
at the singular points hold simultaneously for both  $\vG$ and $\hG$.
\item If $S$ is oriented we require furthermore that also $\hG$ is oriented 
in such a way that $\vG$ and $\hG$ intersect everywhere with intersection number
equal to $1$. 

\end{enumerate}

\begin{lem}\label{sing-f-exist}
  For every ammissible function $\iG$ there are transverse pairs of
  foliations of type $\iG$ on $(S,V)$.
\end{lem}
\Dim It is enough to prove that there exists a vertical foliation of
type $\iG$, for we can take as transverse horizontal foliation the
orhogonal one with respect to a suitable auxiliary riemannian metric
on $S$.  Let $S_Z$ be the surface with boundary obtained by removing
from $S$ a small $2$-disk around every $v\in Z$.  Consider the
foliation on a neighbourhood of $\partial S_Z$ determined by the
$\iG$-local model at singular points.  By a simple variation of Hopf
theorem we realizes that it extends to the whole of $S_Z$ without
introducing new singularities.

\cvd  

\smallskip  

Such transvere pairs of type $\iG$ are considered up homotopy through
transverse pairs of type $\iG$ which is locally an isotopy at the
singular points. We denote by $$\Tt \Pp(S,V,\iG)$$ the so obtained
quotient set. Set
$$ \Tt \Pp(S,V)= \cup_\iG \Tt \Pp(S,V,\iG) \ . $$ Finally denote
by $$\TG(S_V)$$ the quotient set of the set of generic traversing
foliations on $S_V$ considered up to homotopy through traversing (not
necessarily generic) foliations.
\smallskip

\begin{figure}[ht]
\begin{center}
 \includegraphics[width=6cm]{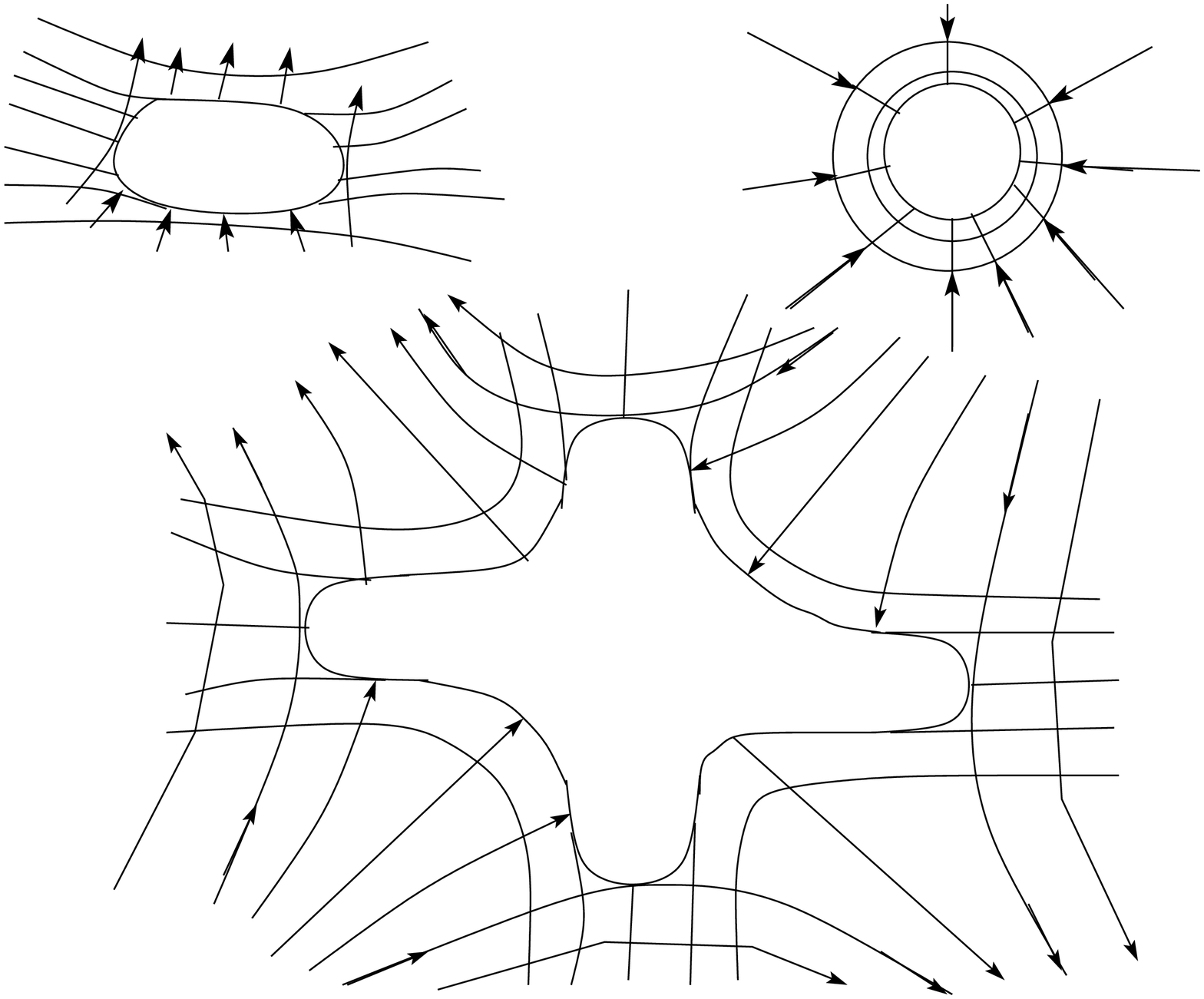}
\caption{\label{SingF} Extendible configurations at $\partial S_V$} 
\end{center}
\end{figure}

The following theorem summarizes some main features   
of the $s$-transit equivalence.
\begin{teo}\label{on-s} For every $(S,V)$:

(1) The correspondence $(T,b)\to \Vv(T,b)$ induces a well defined bijection 
$$\tau: \Ss^{id}(S,V)\to \TG(S_V) \ . $$
\smallskip

(2) For every $(T,b)$ consider the admissible function $\iG_b(v) = 1 -
d_b(v)$.  Then the associated pair of transverse foliations
$(\Vv(T,b),\Hh(T,b))$ on $S_V$ extends to a transverse pair
$(\vG(T,b),\hG(T,b))$ of type $\iG_b$ on $(S,V)$ in such a way that
this induces a well defined map
$$\pG: \Ss^{id}(S,V)\to \Tt \Pp(S,V) \ . $$
\smallskip

(3) Fix a base point $v_0\in V$. Assume that the set of admissible
functions such that $\iG(v_0)=0$ is non empty and denote by
$\Tt\Pp_0(S,V)$ the corresponding subset of $\Tt\Pp(S,V)$.  Then the
set of triangulations $(T,b)$ of $(S,V)$ such that $\iG_b(v_0)=0$ is
non empty and denoting by $\Ss_0^{id}(S,V)$ the corresponding subset
of $\Ss^{id}(S,V)$ we have that the restricted map $\pG:
\Ss_0^{id}(S,V)\to \Tt \Pp_0(S,V)$ is bijective.
\end{teo}
\Dim The fact that the map $\tau$ in (1) is well defined and onto
follows just by looking at the sliding flips and from (3) of
Proposition \ref{traversing}.

Once the extension $(\Vv(T,b),\Hh(T,b)) \to (\vG(T,b),\hG(T,b))$ will
be established just below, the fact that the map $\pG$ of (2) is well
defined follows from the fact that $\tau$ is well defined.

The fact that $\tau$ in (1) is injective as well as item (3) are
actually simpler 2D versions of results established in \cite{BP} and
\cite{BP1} for branched spines of $3$-manifolds with non empty
boundary. In \cite{BP} we essentially considered the case of closed
manifolds, that is when the boundary consists of one $2$-sphere. In
\cite{BP1} we faced the general case with minor changes.  So we limit
here to illustrate the main points, referring to the harder proofs in
3D.
\smallskip

The injectivity of $\tau$ is the 2D analogous of Theorem 4.3.3 of
\cite{BP}.  By transversality we can assume that the homotopy is
generic, that is, it contains only a finite number of non generic
traversing foliations, each one containing one exceptional leaf which
is tangent at two points of $\partial S_V$. Then we analyze how two
generic traversing foliations close to a non generic one are related
to each other and we realizes that the sliding $b$-flips cover all
possible configurations.
\smallskip

As for the extension of $(\Vv,\Hh)$ we look at the configuration of
this foliations at each component $C$ of $\partial S_V$ which is also
the boundary of a small disk $D$ in $S$ centred at one vertex $v\in V$
(see Figure \ref{SingF}).  If $\iG_b(v)=0$, there are exactly two
exceptional leaves of $\Vv$ tangent to $C$; then we can extend $\Vv$
without singularities through $D$ respecting the bicoloring of $C$ and
manage similarly on $\Hh$.  In the other case we easily realize that
up to isotopy the configuration of $(\Vv,\Hh)$ at $C$ is the
restriction of the configuration of the local models that we can
assume to be carried by $D$. The arbitrary choices in implementing the
construction are immaterial if $(\vG,\hG)$ are considered up to kind
of homotopy stated above.
\smallskip

Item (3) is more demanding. First notice that if $\chi(S)=0$ then the
hypothesis can be satisfied also when $|V|=1$, in all other cases
necessarily $|V|\geq 2$. First we prove at the same time that
$\Ss_0^{id}(S,V)$ is non empty and the restricted map $\pG:
\Ss_0^{id}(S,V)\to \Tt \Pp_0(S,V)$ is onto. The key point is to prove
that every vertical foliation $\vG$ occurring in $\Tt \Pp_0(S,V)$ is
realized by means of a triangulation $(T,b)$ with the given
distribution of $d_b(v)$'s.  This is the 2D counterpart of Proposition
5.1.1 of \cite{BP}. The proof is based on Ishii's notion of {\it flow
  spines} \cite{I}.  Let $A:= V \setminus \{v_0\}$ and define $S_A$ as
usual, so that $S_V \subset S_A$. In general $\vG$ is not traversing
$S_A$.  Qualitatively the idea is to detect an embedded $2$-disk $D$
in the interior of $S_A$, centred at some point $v_0'$, such that
$\vG$ becomes traversing $S_{V'}$ and generic, where $V'=\{v_0'\} \cup
A$, and such that only two exceptional leaves are tangent to $\partial
D$. This is carried by a branched triangulation of $(S,V')$. Finally,
via the homogeneity of $S$, we get the desired triangulation $(T,b)$
of $S_V$. The injectivity of the restricted map $\pG$ is the
counterpart of Theorem 5.2.1 of \cite{BP}. The basic idea is to
`cover' any homotopy connecting $\vG(T,b)$ with $\vG(T',b')$ with a
chain of flow-spines connecting $(T,b)$ with $(T',b')$ such that the
traversing foliation associated to one is homotopic through traversing
foliation to the traversing foliations of the subsequent.
 
 \cvd
 
 \begin{remark}\label{infinite-s}{\rm
     Every $\Tt \Pp(S,V,\iG)$ is an affine space over $H_1(S_Z;\Z)$,
     $Z$ being the singular set prescribed by $\iG$. So in general
     $\pG: \Ss_0^{id}(S,V)\to \Tt \Pp_0(S,V)$ is a bijection between
     infinite sets.}
 \end{remark}

  \subsection{On the non ambiguous transit}\label{na} 
 In $3D$ the notion of non ambiguous structure, defined indeed in
 terms of transit of {\it pre-branchings} rather than of branchings \cite{NA}, \cite{3D},
 gives rise to non trivial examples of intrinsic interest.  
 In 2D the intrinsic content of the $na$-relation is not
 so evident.

 \begin{exa}\label{na-ex}{\rm Let $S$ be the torus and $|V|=1$.
 Let $T$ be the triangulation of $(S,V)$ as in the proof of
 Proposition \ref{b-exist}.  One checks by direct inspection that for
 every branching $(T,b)$ no edge of $T$ supports a non ambiguous flip.
 This holds for {\it every} triangulation $T'$ of $(S,V)$ because all
 these triangulations are equivalent to each other up to
 diffeomorphism of $(S,V)$.  Hence in this case the $na$-transit
 equivalence is nothing else than the identity relation. On the other
 hand , we check that the branchings on
 $T$ which share the same decomposition $S=S_+\cup S_-$ are not
 $s$-equivalent to each other.}
  \end{exa}
 
 Referring to Proposition \ref{S+}, one would conjecture, better ask
 whether two branched triangulations $(T,b)$ and $(T',b')$ of $(S,V)$
 are $na$-equivalent if and only if they are $s$-equivalent and share
 the decomposition $S=S_+\cup S_-$, provided that $S$ is oriented.
 However, while by (3) of Theorem \ref{on-s} the quotient set
 $\Ss^{id}_0(S,V)$ has a nice intrinsic content, the
 decomposition $S=S_+\cup S_-$ is not very transparent.
 Nevertheless, again when $S$ is oriented, we will point out a
 further structure preserved
 by the $na$-transits with a bit more intrinsic flavour.
We set
 $H^\Delta_*$ to denote the simplicial (or cellular) homology of a
 complex, $H_*$ the singular homology of the underlying topological
 space. Similarly we set $H^*_\Delta$, $H^*$ for the cohomology.  The
 inclusion of $\theta$ in $S_V$ induces an isomorphism
$$H^\Delta_1(\theta,b;\R) \cong H_1(S_V;\R)$$
and via elementary Poincar\'e duality we have
$$ H^\Delta_1(\theta,b;\R)\cong H^1_\Delta(T,b;\R)\cong H^1(S;\R)\ . $$
Moreover, $H^\Delta_1(\theta,b;\R)=Z^\Delta_1(\theta,b;\R)$
where this last denotes the space of $1$-cycles on $(\theta,b)$.
Every $z\in Z^\Delta_1(\theta,b;\R)$ consists in giving each $b$-oriented edge
$e$ of $(\theta,b)$ a weight $z(e)\in \R$ in such a way that at every
switching vertex of $\theta$ the three weights around the vertex verify
the corresponding switching condition of the form $z(e_0)=z(e_1)+z(e_2)$.
These cycles transit along every $b$-flip, so that for every composite
$b$-transit $(T,b)\Rightarrow (T',b')$ it is defined an isomorphism
$$\alpha: Z^\Delta_1(\theta,b;\R)\to Z^\Delta_1(\theta',b';\R) \ . $$
Set
$$\Mm=\Mm_{(T,b)} = \{z\in Z^\Delta_1(\theta,b;\R)| \ \forall 
e\in \theta^{(1)}, \ z(e)\geq 0\}; \
\Mm^+ = \{z\in Z^\Delta_1(\theta,b;\R)| \ \forall 
e\in \theta^{(1)}, \ z(e)> 0\}\  . $$
Every $z\in \Mm$ can be interpreted as a {\it transverse measure}
on the horizontal foliation $\Hh$ or on the singular foliation $\hG$.
By taking into account the arbitrary choices in the realizations of 
$(\Vv$, $\Hh)$ and $(\vG,\hG)$
we radily have

\begin{prop}\label{same measure}
  (1) For every $z\in \Mm$, the measured foliations $(\Hh,z)$,
  $(\hG,z)$ are uniquely detemined up to {\rm measure equivalence} (in
  particular this means that $(\hG,z)$ well define a {\rm measure
    spectrum} on the set $\Ss$ of isotopy classes of simple closed
  curves on $S$).
  
  (2) If we denote by $\Mm(\hG)$ the set of transverse measures 
on $\hG$ up to measure equivalence, the above correspondence 
well defines a map
$$ \mG=\mG_{(T,b)}: \Mm_{(T,b)}  \to \Mm(\hG) . \ $$
\end{prop}

\cvd

After a look at the $na$-transits we readily have

\begin{prop}
  If $(T,b)$ and $(T',b')$ are $na$-transit equivalent, then the maps
  $\mG_{(T,b)}$ and $\mG_{(T',b')}$ have the same image. More
  precisely, there is a bijection $\alpha: \Mm_{(T,b)}\to
  \Mm_{(T',b')}$ such that $\alpha(\Mm^+_{(T,b)})= \Mm^+_{(T',b')}$
  and $\mG_{(T,b)}= \mG_{(T',b')}\circ \alpha$.
\end{prop}

\cvd

\end{document}